\documentclass[11pt]{article}
\usepackage{latexsym,amscd,amsmath,amssymb}

\begin{document}
\newcommand{\Cn}{\mathbb{C}^n}
\newcommand{\he}[2]{\sum\limits_{#1}^{#2}}
\newcommand{\lj}[2]{\prod\limits_{#1}^{#2}}
\newcommand{\mo}[1]{\parallel #1\parallel_p}
\newcommand{\jz}[5]{\left( \begin{array}{c}#1 \\ #2 \\#3\\#4\\#5\end{array} \right)}
\newcommand{\jw}[4]{\left( \begin{array}{c}#1 \\ #2 \\#3\\#4\end{array} \right)}
\newcommand{\jwz}[2]{\left( \begin{array}{c}#1 \\ #2 \end{array} \right)}
\newcommand{\lwjz}[5]{\left( \begin{array}{lllll}#1 \\ #2 \\#3\\#4\\#5\end{array} \right)}
\newcommand{\hwjz}[5]{\begin{array}{ccccc}#1 & #2 &#3&#4&#5\end{array}}
\newcommand{\hsjz}[4]{\begin{array}{cccc}#1 & #2 &#3&#4\end{array}}
\newcommand{\jd}[2]{\mid #1\mid }
\newcommand{\fd}[4]{\left \{ \begin{array}{ll}
#1&#2\\#3&#4\end{array}\right.}

\title{\bf Criteria of Biholomorphic Convex Mappings on the bounded convex balanced domain $D_{p}^n$ }

\author{Ni Li and Ming-Sheng Liu\thanks{Correspondence should be addressed to Ming-Sheng Liu;\, liumsh@scnu.edu.cn}\\
 {\it School of Mathematical Sciences, South China Normal University,} \\
 {\it Guangzhou 510631, Guangdong, People's Republic of China}}

\date{}
\maketitle
\pagestyle{plain}
\renewcommand{\thefootnote}{\fnsymbol{footnote}}

\footnote[0]{This research is supported by Guangdong Natural Science Foundation (Grant No.2014A030307016, 2014A030313422).}

\begin{center}
\begin{minipage}{13.5cm}
{\bf Abstract.}\quad In this paper, we first establish several general sufficient conditions for the biholomorphic convex mappings on the bounded convex balanced domain $D_{p}^n(p_{j}\geq 2,j=1,\cdots,n)$ in $C^{n}$, which extend some related  results of earlier authors. From these, some concrete examples of biholomorphic convex mappings on $D_{p}^n$ are also provided.
\end{minipage}
\end{center}

\section{Introduction and Preliminaries}
\quad\, Suppose that $C^{n}$ is the vector space of $n$ complex variables $z=(z_{1},z_{2},\cdots,z_{n})$ with the Euclidean inner product $\langle z,w\rangle =\he{j=1}{n}z_{j}\bar{w}_{j}$, where $w=(w_{1},w_{2},\cdots,w_{n})\in C^{n}$. If $\Omega$ is a domain in $C^n$, and for every $z\in \Omega$, $|\lambda|\leq 1$, we have $\lambda z\in \Omega$, then we call $\Omega$ a balanced domain. Its Minkowski functional is defined by
$$
\rho(z)=\inf \{t>0,\frac{z}{t}\in \Omega\},\quad z\in C^{n}.
$$

Assume $\Omega$ is a bounded convex balanced domain in $C^{n}$, let $\rho(z)$ be its Minkowski functional, then $\rho(\cdot)$ is a norm of $C^{n}$, $\rho(\lambda z)=|\lambda|\rho(z)$ for every $\lambda\in C, z\in C^{n}$, and
$$
\Omega=\{z\in C^{n}:\rho(z)< 1\},\quad
$$
and $\rho(z)=0$ if and only if $z=0$.

Let $p=(p_1, \cdots, p_n)$ with $p_j> 1\, (j=1,2,\cdots,n)$, and
\begin{eqnarray}
D_{p}^n=\{(z_{1},z_{2},\cdots,z_{n})\in C^{n}:\he{j=1}{n}|z_{j}|^{p_{j}}< 1\},
\label{liu1}
\end{eqnarray}
then $D_{p}^n$ is a bounded convex balanced domain in $C^{n}$, and its Minkowski functional $\rho(z)$ satisfies the following equality
\begin{eqnarray}
\he{j=1}{n}\bigg|\frac{z_{j}}{\rho(z)}\bigg|^{p_{j}}=1.
\label{liu2}
\end{eqnarray}

When $p_1=\cdots=p_n=p$, we denote $D_{p}^n$ by $B^{n}_{p}$, at this time, we have $\rho(z)=\sqrt[p]{|z_{1}|^{p}+\cdots+|z_{n}|^{p}}$. In particular, let $U=B^{1}_{p}$ be the unit disk in the complex plane $C$.

Let $\Omega$ be a domain in $C^{n}$, a mapping $f: \Omega\rightarrow C^{n}$ is said to be locally biholomorphic in $\Omega$ if $f$ has a locally inverse at each point $z\in \Omega$ or, equivalently, if the first Fr\'{e}chet derivative
$$
Df(z)=\bigg(\frac{\partial f_{j}(z)}{\partial z_{k}}\bigg)_{1 \leq j,k\leq n}
$$
is nonsingular at each point in $\Omega$. 
The second Fr\'{e}chet derivative of a mapping $f:\Omega\rightarrow C^{n}$ is a symmetric bilinear operator $D^{2}f(z)(\cdot,\cdot)$ on $C^{n}\times C^{n}$, and $ D^{2}f(z)(z,\cdot)$ is the linear operator
obtained by restricting $D^2f(z)$ to $\{z\}\times C^n$. The matrix
representation of $D^2f(z)(b,\cdot)$ is
$$
D^{2}f(z)(b,\cdot)=\bigg(\he{l=1}{n}\frac{\partial^{2}f_{j}(z)}{\partial z_{k}\partial z_{l}}b_{l}\bigg)_{1\leq j,k\leq n}\\
$$
where $f(z)=(f_{1}(z),\cdots,f_{n}(z)),b=(b_{1},\cdots,b_{n})\in C^{n}$.

Let $N(D_{p}^n)$ be the class of all locally biholomorphic mappings $f(z)=(f_{1}(z),\cdots,f_{n}(z)): D_{p}^n\rightarrow C^{n}$ such that $f(0)=0,Df(0)=I$, where $z=(z_{1},\cdots,z_{n})\in C^{n} $ and $I$ is the unit matrix of $ n\times n$. If $f\in N(D_{p}^n)$ is a biholomorphic mapping on $D_{p}^n$ and $f(D_{p}^n)$ is a convex domain in $C^{n}$, then we say that $f$ is a biholomorphic convex mappings on $D_{p}^n$. Let $K(D_{p}^n)$ denote the class of all biholomorphic convex mappings
on $D_p^n$ with $f(0)=0, Df(0)=I$.

It is not easy to construct concrete biholomorphic convex
mappings on some domains in $C^{n}$, even on the unit ball $B^{n}_{2}$.
In 1995, Roper and Suffridge\cite{rs1} proved that: ~If $f\in K$ and $F(z)=(f(z_{1}),
\sqrt{f'(z_{1})}z_{0})$, where $z=(z_1,z_0)\in B_2^n, z_1\in U, z_0=(z_2,\cdots, z_n) \in C^{n-1}$, then $F\in K(B^{n}_{2})$.
Which is popularly referred to as the Roper-Suffridge operator. Using this operator, we may construct a lot of concrete biholomorphic
convex mappings on $B_n^2$. Gong and Liu \cite{gl} generalized the Roper-Suffridge operator to the Reinhardt domain $D_{p}=\{z=(z_{1},z_{2},\cdots,z_{n})\in C^{n}: \he{j=1}{n}|z_{j}|^{p_{j}}< 1\}$, where $p_{1}=2,p_{2}\geq 1, p_{3}\geq 1, \cdots, p_{n}\geq 1$. However, according to the result in \cite{g, lz}, none of these concrete examples
belongs to $K(D_p^n)(p_j>2, j=1, 2, \cdots, n)$. In 1999, Roper and Suffridge\cite{rs2} gave two concrete examples of biholomorphic convex
mappings on $B^{2}_{p}$ in $C^{n}$,  Liu and Zhu\cite{lz2006} gave some sufficient conditions for biholomorphic convex mappings on $B^{n}_{p}$.  Hamada and Kohr\cite{gk0}, Zhu\cite{z1} gave a necessary and sufficient
condition for biholomorphic convex mappings on bounded convex balanced domain $D_{p}^n$ as follows.

{\bf Theorem A}\cite{gk0,z1}\quad Suppose that $ p_{j}\geq 2 (j=1,2,...,n)$, $\rho(z)$ is the Minkowski functional of $D_{p}^n$, and $ f\in N(D_{p}^n)$. Then $f\in K(D_{p}^n)$ if and only if for any $z=(z_{1},z_{2},\cdots,z_{n})\in D_{p}^n\backslash\{0\}$, and $b=(b_{1},b_{2},\cdots,b_{n})\in C^{n}\backslash{\{0\}}$ such that
$$
\mbox{Re}\bigg\{\he{j=1}{n}p_{j}\bigg|\frac{z_{j}}{\rho(z)}\bigg|^{p_{j}}\frac{b_{j}}{z_{j}}\bigg\}=0,
$$
we have
\begin{eqnarray*}
J_{f}(z,b)&=&\mbox{Re}\bigg\{\he{j=1}{n}\frac{p^{2}_{j}}{2}\frac{|z_{j}|^{p_{j}-2}}{\rho(z)^{p_{j}}}|b_{j}|^{2}+\he{j=1}{n}p_{j}(\frac{p_{j}}{2}-1)
\Big|\frac{z_{j}}{\rho(z)}\Big|^{p_{j}}(\frac{b_{j}}{z_{j}})^{2}\\
&&-2\he{j=1}{n}\frac{p_{j}}{\rho(z)}\Big|\frac{z_{j}}{\rho(z)}\Big|^{p_{j}}\bigg\langle
Df(z)^{-1}D^{2}f(z)(b,b),\frac{\partial \rho}{\partial\overline{z}}\bigg\rangle
\bigg\}\geq 0.
\end{eqnarray*}

Liu and Zhu\cite{lz2009} had established some sufficient conditions
of biholomorphic convex mappings on $D_{p}^n$ for mappings of the following forms
$$
f(z)=(f_1(z_{1})+g_{1}(z_{n}), f_2(z_{2})+g_{2}(z_{n}),\cdots , f_{n-1}(z_{n-1})+g_{n-1}(z_{n}), f_n(z_{n})),
$$
and
$$
f(z)=(f_1(z_1,z_2,\cdots, z_n), f_2(z_2),\cdots, f_n(z_n)).
$$

Liu and Li \cite{ll2014} did the further promotion on $D_{p}^n$, they had established a sufficient condition
of biholomorphic convex mappings on $D_{p}^n$ for mappings of the following form
$$
f(z)=(f_1(z_1,z_2,\cdots, z_n),  f_2(z_{2})+g_{2}(z_{n}),\cdots , f_{n-1}(z_{n-1})+g_{n-1}(z_{n}), f_n(z_{n})).
$$

Now we pose two problems as follows:

Problem I: Can we get some sufficient conditions such that mapping of the form :
$$
f(z)=(f_{1}(z_{1},z_{2},\cdots,z_{n}),f_{2}(z_{2},\cdots,z_{n}),\cdots,f_{n-1}(z_{n-1},z_{n}),f_{n}(z_{n}))
$$
is a biholomorphic convex mapping on $D_{p}^n$?

Problem II: Can we get some sufficient conditions such that the mapping
$$
f(z)=(f_{1}(z_{1},z_{2}),f_{2}(z_{1},z_{2})),\quad z=(z_1,z_2)
$$
is a biholomorphic convex mapping on $D_{p}^2$ in $C^2$?

The object of this paper is to give partial answers to the above problems. From these, we may construct some concrete biholomorphic convex mappings on $D_{p}^n$.

\section{Main results}
\quad\, In order to give a partial answer to Problem I, we first establish the following theorem.

{\bf Theorem 1}\quad Suppose that $n\geq 2,  p_j\ge 2,j=1,2,\cdots,n$, and $f_{j}:U\rightarrow C$ is holomorphic with $f_{j}(0)=0,f'_{j}(0)=1(j=3,\cdots,n-1,n)$, $f_{1}(z_{1},\cdots,z_{n}):D_{p}^n\rightarrow C$ is holomorphic with $f_{1}(0,0,\cdots,0)=0,\frac{\partial f_{1}}{\partial z_{1}}(0,0,\cdots,0)=1,\frac{\partial f_{1}}{\partial z_{l}}(0,0,\cdots,0)=0(l=2,3,\cdots,n)$. Let
$$
f(z)=(f_{1}(z_{1},z_{2},\cdots,z_{n}),f_{2}(z_{2},z_{3},\cdots,z_{n}),
f_{3}(z_{3}),\cdots,f_{n-1}(z_{n-1}),f_{n}(z_{n})).
$$
If $f(z)$ satisfies the following conditions
\begin{eqnarray*}
&(1)&\frac{\partial f_{1}}{\partial z_{1}}\frac{\partial f_{2}}{\partial z_{2}}\prod_{j=3}^{n}f'_{j}(z_{j})\neq 0,\quad |z_{j}f''_{j}(z_{j})|\leq|f'_{j}(z_{j})|(j=3,4,\cdots,n);\\
&(2)&\he{l=1}{n}|z_{1}\frac{\partial^{2}f_{1}}{\partial z_{1}\partial z_{l}}|\leq|\frac{\partial f_{1}}{\partial z_{1}}|,\quad \he{l=2}{n}|z_{2}\frac{\partial^{2}f_{2}}{\partial z_{2}\partial z_{l}}|\leq|\frac{\partial f_{2}}{\partial z_{2}}|;\\
&(3)&p_{1}\bigg(\he{l=1}{n}\bigg|\frac{\frac{\partial^{2}f_{1}}{\partial z_{2}\partial z_{l}}}{\frac{\partial f_{1}}{\partial z_{1}}}\bigg|+\he{l=2}{n}\bigg|\frac{\frac{\partial f_{1}}{\partial z_{2}}}{\frac{\partial f_{2}}{\partial z_{2}}}\bigg|\bigg|\frac{\frac{\partial^{2}f_{2}}{\partial z_{2}\partial z_{l}}}{\frac{\partial f_{1}}{\partial z_{1}}}\bigg|\bigg)\leq\bigg(1-\he{l=2}{n}\bigg|\frac{\frac{z_{2}\partial^{2}f_{2}}{\partial z_{2}\partial z_{l}}}{\frac{\partial f_{2}}{\partial z_{2}}}\bigg|\bigg)p_{2}|z_{2}|^{p_{2}-2}; \\
&(4)&p_{1}\bigg(\he{l=1}{n}\bigg|\frac{\frac{\partial^{2}f_{1}}{\partial z_{j}\partial z_{l}}}{\frac{\partial f_{1}}{\partial z_{1}}}\bigg|+\he{l=2}{n}\bigg|\frac{\frac{\partial f_{1}}{\partial z_{2}}}{\frac{\partial f_{2}}{\partial z_{2}}}\bigg|\bigg|\frac{\frac{\partial^{2}f_{2}}{\partial z_{j}\partial z_{l}}}{\frac{\partial f_{1}}{\partial z_{1}}}\bigg|+\bigg|\frac{\frac{\partial f_{1}}{\partial z_{2}}}{\frac{\partial f_{1}}{\partial z_{1}}}\bigg|\bigg|\frac{\frac{\partial f_{2}}{\partial z_{j}}}{\frac{\partial f_{2}}{\partial z_{2}}}\bigg|\bigg|\frac{f''_{j}(z_{j})}{f'_{j}(z_{j})}\bigg|+\bigg|\frac{\frac{\partial f_{1}}{\partial z_{j}}}{\frac{\partial f_{1}}{\partial z_{1}}}\bigg|\bigg|\frac{f''_{j}(z_{j})}{f'_{j}(z_{j})}\bigg|\bigg)\\
&&+p_{2}\bigg(\he{l=2}{n}\bigg|\frac{\frac{\partial^{2}f_{2}}{\partial z_{j}\partial z_{l}}}{\frac{\partial f_{2}}{\partial z_{2}}}\bigg|+\bigg|\frac{\frac{\partial f_{2}}{\partial z_{j}}}{\frac{\partial f_{2}}{\partial z_{2}}}\bigg|\bigg|\frac{f''_{j}(z_{j})}{f'_{j}(z_{j})}\bigg|\bigg)\leq p_{j}|z_{j}|^{p_{j}-2}\bigg(1-\bigg|z_{j}\frac{f''_{j}(z_{j})}{f'_{j}(z_{j})}\bigg|\bigg) (j=3,4,\cdots,n)
\end{eqnarray*}
for all $z=(z_1,\cdots, z_n)\in D_p^n\setminus\{0\}$, then $f\in K(D_p^n)$.

{\bf Proof}\quad By direct computation of the Fr\'{e}chet derivatives of $f(z)$ , we get that
\begin{gather*}
Df(z)=\begin{pmatrix}\frac{\partial f_{1}}{\partial z_{1}}&\frac{\partial f_{1}}{\partial z_{2}}&\cdots&\frac{\partial f_{1}}{\partial z_{n-1}}&\frac{\partial f_{1}}{\partial z_{n}}\\
0&\frac{\partial f_{2}}{\partial z_{2}}&\cdots&\frac{\partial f_{2}}{\partial z_{n-1}}&
\frac{\partial f_{2}}{\partial z_{n}}\\
\cdots&\cdots&\cdots&\cdots&\cdots\\
0&0&\cdots&f'_{n-1}(z_{n-1})&0\\
0&0&\cdots&0&f'_{n}(z_{n})\end{pmatrix}
\end{gather*}
\begin{gather*}
Df(z)^{-1}=\begin{pmatrix}\frac {1}{\frac{\partial f_{1}}{\partial z_{1}}}&-\frac{\frac{\partial f_{1}}{\partial z_{2}}}{\frac{\partial f_{1}}{\partial z_{1}}\frac{\partial f_{2}}{\partial z_{2}}}&\frac{\frac{\partial f_{1}}{\partial z_{2}}\frac{\partial f_{2}}{\partial z_{3}}-\frac{\partial f_{1}}{\partial z_{3}}\frac{\partial f_{2}}{\partial z_{2}}}{\frac{\partial f_{1}}{\partial z_{1}}\frac{\partial f_{2}}{\partial z_{2}}f'_{3}(z_{3})}&\cdots&\frac{\frac{\partial f_{1}}{\partial z_{2}}\frac{\partial f_{2}}{\partial z_{n-1}}-\frac{\partial f_{1}}{\partial z_{n-1}}\frac{\partial f_{2}}{\partial z_{2}}}{\frac{\partial f_{1}}{\partial z_{1}}\frac{\partial f_{2}}{\partial z_{2}}f'_{n-1}(z_{n-1})}& \frac{\frac{\partial f_{1}}{\partial z_{2}}\frac{\partial f_{2}}{\partial z_{n}}-\frac{\partial f_{1}}{\partial z_{n}}\frac{\partial f_{2}}{\partial z_{2}}}{\frac{\partial f_{1}}{\partial z_{1}}\frac{\partial f_{2}}{\partial z_{2}}f'_{n}(z_{n})}\\
0&\frac {1}{\frac{\partial f_{2}}{\partial z_{2}}}&-\frac{\frac{\partial f_{2}}{\partial z_{3}}}{\frac{\partial f_{2}}{\partial z_{2}}f'_{3}(z_{3})}&\cdots&-\frac{\frac{\partial f_{2}}{\partial z_{n-1}}}{\frac{\partial f_{2}}{\partial z_{2}}f'_{n-1}(z_{n-1})}&-\frac{\frac{\partial f_{2}}{\partial z_{n}}}{\frac{\partial f_{2}}{\partial z_{2}}f'_{n}(z_{n})}\\
\cdots&\cdots&\cdots&\cdots&\cdots&\cdots\\
0&0&\cdots&\cdots&\frac{1}{f'_{n-1}(z_{n-1})}&0 \\
0&0&\cdots&0&0&\frac {1}{f'_{n}(z_{n})}\end{pmatrix}
\end{gather*}
\begin{eqnarray*}
D^{2}f(z)(b,b)=\begin{pmatrix}
\he{l=1}{n}\frac{\partial^{2}f_{1}}{\partial z_{1}\partial z_{l}}b_{l}&\he{l=1}{n}\frac{\partial^{2}f_{1}}{\partial z_{2}\partial z_{l}}b_{l}&\cdots&\he{l=1}{n}\frac{\partial^{2}f_{1}}{\partial z_{n}\partial z_{l}}b_{l}\\
0&\he{l=2}{n}\frac{\partial^{2}f_{2}}{\partial z_{2}\partial z_{l}}b_{l}&\cdots&\he{l=2}{n}\frac{\partial^{2}f_{2}}{\partial z_{n}\partial z_{l}}b_{l}\\
\cdots&\cdots&\cdots&\cdots\\
0&0&\cdots&0\\
0&0&\cdots&f''_{n}(z_{n})b_{n}
\end{pmatrix}
\begin{pmatrix}
b_{1}\\b_{2}\\\cdots\\b_{n-1}\\b_{n}
\end{pmatrix}
=\begin{pmatrix}\he{j=1}{n}\he{l=1}{n}\frac{\partial^{2}f_{1}}{\partial z_{j}\partial z_{l}}b_{l}b_{j}\\
\he{j=2}{n}\he{l=2}{n}\frac{\partial^{2}f_{2}}{\partial z_{j}\partial z_{l}}b_{l}b_{j}\\
f''_{3}(z_{3})b^{2}_{3}\\
\cdots\\
f''_{n}(z_{n})b^{2}_{n}\end{pmatrix}.
\end{eqnarray*}

From (\ref{liu2}) and $|z_j|^{p_j}={z_j}^\frac{p_j}{2}\,\,{\overline{z_j}}^\frac{p_j}{2}$, direct computation yields
\begin{equation}
\frac{\partial\rho}{\partial\overline{z_l}}
=\frac{p_l|z_l|^{p_l}}{2\overline{z_l}\rho(z)^{p_l-1}\he{j=1}{n}p_j\big|\frac{z_j}{\rho(z)}\big|^{p_j}},
\label{liu3}
\end{equation}
and $\Big|\frac{z_j}{\rho(z)}\Big|\leq1(j=1, 2, \cdots, n)$ for all
$z=(z_1,z_2,\cdots,z_n)\in D_p^n$.

Taking $z=(z_{1}, \cdots ,z_{n})\in D^n_p\setminus\{0\}, b=(b_{1},\cdots, b_{n})\in C^n\setminus\{0\}$ such that
$$\mbox{Re}\bigg\{\he{j=1}{n}p_{j}\bigg|\frac{z_{j}}{\rho(z)}\bigg|^{p_{j}}\frac{b_{j}}{z_{j}}\bigg\}=0,$$
by applying the fact $0<\rho(z)<1$ and the hypothesis of Theorem 1, we have
\begin{eqnarray*}
&J_{f}(z,b)&\geq\he{j=1}{n}\frac{p_{j}|z_{j}|^{p_{j}-2}}{\rho(z)^{p_{j}}}|b_{j}|^{2}-2\he{j=1}{n}\frac{p_{j}}{\rho(z)}\Big|\frac{z_{j}}{\rho(z)}\Big|^{p_{j}}
\mbox{Re}\langle Df(z)^{-1}D^{2}f(z)(b,b), \frac{\partial \rho}{\partial\overline{z}}\rangle\\
&=&\he{j=1}{n}\frac{p_{j}|z_{j}|^{p_{j}-2}}{\rho(z)^{p_{j}}}|b_{j}|^{2}-\mbox{Re}\bigg\{\frac{1}{\frac{\partial f_{1}}{\partial z_{1}}}\bigg[\he{j=1}{n}\he{l=1}{n}\frac{\partial^{2}f_{1}}{\partial z_{j}\partial z_{l}}b_{l}b_{j}-\he{j=2}{n}\he{l=2}{n}\frac{\frac{\partial f_{1}}{\partial z_{2}}}{\frac{\partial f_{2}}{\partial z_{2}}}\frac{\partial^{2}f_{2}}{\partial z_{j}\partial z_{l}}b_{j}b_{l}\\
&&+\he{j=3}{n}\bigg(\frac{\frac{\partial f_{1}}{\partial z_{2}}}{\frac{\partial f_{2}}{\partial z_{2}}}\frac{\partial f_{2}}{\partial z_{j}}-\frac{\partial f_{1}}{\partial z_{j}}\bigg)\frac{f''_{j}(z_{j})}{f'_{j}(z_{j})}b^{2}_{j}\bigg]\frac{p_{1}|z_{1}|^{p_{1}}}{z_{1}\rho(z)^{p_{1}}}\\
&&+\bigg
[\he{j=2}{n}\he{l=2}{n}\frac{\frac{\partial^{2}f_{2}}{\partial z_{j}\partial z_{l}}}{\frac{\partial f_{2}}{\partial z_{2}}}b_{l}b_{j}-\he{j=3}{n}\frac{\frac{\partial f_{2}}{\partial z_{j}}}{\frac{\partial f_{2}}{\partial z_{2}}}\frac{f''_{j}(z_{j})}{f'_{j}(z_{j})}b^{2}_{j}\bigg]\frac{p_{2}|z_{2}|^{p_{2}}}{z_{2}\rho(z)^{p_{2}}}
+\he{j=3}{n}\frac{f''_{j}(z_{j})}{f'_{j}(z_{j})}b^{2}_{j}\frac{p_{j}|z_{j}|^{p_{j}}}{z_{j}\rho(z)^{p_{j}}}\bigg\}\\
&\geq&\he{j=1}{n}\frac{p_{j}|z_{j}|^{p_{j}-2}}{\rho(z)^{p_{j}}}|b_{j}|^{2}-\frac{1}{|\frac{\partial f_{1}}{\partial z_{1}}|}\bigg[\he{j=1}{n}\he{l=1}{n}\bigg|\frac{\partial^{2}f_{1}}{\partial z_{j}\partial z_{l}}\bigg||b_{j}|^{2}+\he{j=2}{n}\he{l=2}{n}\bigg|\frac{\frac{\partial f_{1}}{\partial z_{2}}}{\frac{\partial f_{2}}{\partial z_{2}}}\bigg|\bigg|\frac{\partial^{2}f_{2}}{\partial z_{j}\partial z_{l}}\bigg||b_{j}|^{2}\\
&&+\he{j=3}{n}\bigg(\frac{|\frac{\partial f_{1}}{\partial z_{2}}||\frac{\partial f_{2}}{\partial z_{j}}|}{|\frac{\partial f_{2}}{\partial z_{2}}|}+\bigg|\frac{\partial f_{1}}{\partial z_{j}}\bigg|\bigg)\bigg|\frac{f''_{j}(z_{j})}{f'_{j}(z_{j})}\bigg||b_{j}|^{2}\bigg]\frac{p_{1}|z_{1}|^{p_{1}-1}}{\rho(z)^{p_{1}}}-\he{j=3}{n}\bigg|\frac{f''_{j}(z_{j})}{f'_{j}(z_{j})}\bigg|
|b_{j}|^{2}\frac{p_{j}|z_{j}|^{p_{j}-1}}{\rho(z)^{p_{j}}}\\
&&-\bigg[\he{j=2}{n}\he{l=2}{n}\bigg|\frac{\frac{\partial^{2}f_{2}}{\partial z_{j}\partial z_{l}}}{\frac{\partial f_{2}}{\partial z_{2}}}\bigg||b_{j}|^{2}+\he{j=3}{n}\bigg|\frac{\frac{\partial f_{2}}{\partial z_{j}}}{\frac{\partial f_{2}}{\partial z_{2}}}\bigg|\bigg|\frac{f''_{j}(z_{j})}{f'_{j}(z_{j})}\bigg||b_{j}|^{2}
\bigg]\frac{p_{2}|z_{2}|^{p_{2}-1}}{\rho(z)^{p_{2}}}\\
&=&|b_{1}|^{2}\frac{p_{1}|z_{1}|^{p_{1}-2}}{\rho(z)^{p_{1}}}\bigg(1-\he{l=1}{n}
\bigg|\frac{\frac{z_{1}\partial^{2}f_{1}}{\partial z_{1}\partial z_{l}}}{\frac{\partial f_{1}}{\partial z_{1}}}\bigg|\bigg)+\frac{|b_{2}|^{2}}{\rho(z)}\bigg[\bigg(1-\he{l=2}{n}\bigg|\frac{\frac{z_{2}\partial^{2}f_{2}}{\partial z_{2}\partial z_{l}}}{\frac{\partial f_{2}}{\partial z_{2}}}\bigg|\bigg)\frac{p_{2}|z_{2}|^{p_{2}-2}}{\rho(z)^{p_{2}-1}}\\
&&-p_{1}\Big(\frac{|z_{1}|}{\rho(z)}\Big)^{p_{1}-1}
\bigg(\he{l=1}{n}\bigg|\frac{\frac{\partial^{2}f_{1}}{\partial z_{2}\partial z_{l}}}{\frac{\partial f_{1}}{\partial z_{1}}}\bigg|+\he{l=2}{n}\bigg|\frac{\frac{\partial f_{1}}{\partial z_{2}}}{\frac{\partial f_{2}}{\partial z_{2}}}\bigg|
\bigg|\frac{\frac{\partial^{2}f_{2}}{\partial z_{2}\partial z_{l}}}{\frac{\partial f_{1}}{\partial z_{1}}}\bigg|
\bigg)\bigg]+\he{j=3}{n}\frac{|b_{j}|^{2}}{\rho(z)}\bigg[\frac{p_{j}|z_{j}|^{p_{j}-2}}{\rho(z)^{p_{j}-1}}
\bigg(1-\bigg|\frac{z_{j}f''_{j}(z_{j})}{f'_{j}(z_{j})}\bigg|\bigg)\\
&&-p_{1}\Big(\frac{|z_{1}|}{\rho(z)}\Big)^{p_{1}-1}
\bigg(\he{l=1}{n}\bigg|\frac{\frac{\partial^{2}f_{1}}{\partial z_{j}\partial z_{l}}}{\frac{\partial f_{1}}{\partial z_{1}}}\bigg|+\he{l=2}{n}\bigg|\frac{\frac{\partial f_{1}}{\partial z_{2}}}{\frac{\partial f_{2}}{\partial z_{2}}}\bigg|\bigg|\frac{\frac{\partial^{2}f_{2}}{\partial z_{j}\partial z_{l}}}{\frac{\partial f_{1}}{\partial z_{1}}}\bigg|+\bigg|\frac{\frac{\partial f_{1}}{\partial z_{2}}}{\frac{\partial f_{1}}{\partial z_{1}}}\bigg|
\bigg|\frac{\frac{\partial f_{2}}{\partial z_{j}}}{\frac{\partial f_{2}}{\partial z_{2}}}\bigg|
\bigg|\frac{f''_{j}(z_{j})}{f'_{j}(z_{j})}\bigg|+\bigg|\frac{\frac{\partial f_{1}}{\partial z_{j}}}{\frac{\partial f_{1}}{\partial z_{1}}}\bigg|\bigg|\frac{f''_{j}(z_{j})}{f'_{j}(z_{j})}\bigg|\bigg)\\
&&-p_{2}\Big(\frac{|z_{2}|}{\rho(z)}\Big)^{p_{2}-1}
\bigg(\he{l=2}{n}\bigg|\frac{\frac{\partial^{2}f_{2}}{\partial z_{j}\partial z_{l}}}{\frac{\partial f_{2}}{\partial z_{2}}}\bigg|+\bigg|\frac{\frac{\partial f_{2}}{\partial z_{j}}}{\frac{\partial f_{2}}{\partial z_{2}}}\bigg|\bigg|\frac{f''_{j}(z_{j})}{f'_{j}(z_{j})}\bigg|\bigg)\bigg]\\
&\geq&|b_{1}|^{2}\frac{p_{1}|z_{1}|^{p_{1}-2}}{\rho(z)^{p_{1}}}\bigg(1-\he{l=1}{n}
\bigg|\frac{z_{1}\frac{\partial^{2}f_{1}}{\partial z_{1}\partial z_{l}}}{\frac{\partial f_{1}}{\partial z_{1}}}\bigg|\bigg)+|b_{2}|^{2}\bigg[p_{2}|z_{2}|^{p_{2}-2}\bigg(1-\he{l=2}{n}\bigg
|\frac{z_{2}\frac{\partial^{2}f_{2}}{\partial z_{2}\partial z_{l}}}{\frac{\partial f_{2}}{\partial z_{2}}}\bigg|
\bigg)\\
&&-p_{1}\bigg(\he{l=1}{n}\bigg|\frac{\frac{\partial^{2}f_{1}}{\partial z_{2}\partial z_{l}}}{\frac{\partial f_{1}}{\partial z_{1}}}\bigg|+\he{l=2}{n}\bigg|\frac{\frac{\partial f_{1}}{\partial z_{2}}}{\frac{\partial f_{2}}{\partial z_{2}}}\bigg|\bigg|\frac{\frac{\partial^{2}f_{2}}{\partial z_{2}\partial z_{l}}}{\frac{\partial f_{1}}{\partial z_{1}}}\bigg|\bigg)\bigg]+\he{j=3}{n}|b_{j}|^{2}\bigg[p_{j}|z_{j}|^{p_{j}-2}\bigg
(1-\bigg|\frac{z_{j}f''_{j}(z_{j})}{f'_{j}(z_{j})}\bigg|\bigg)\\
\end{eqnarray*}
\begin{eqnarray*}
&&-p_{1}\bigg(\he{l=1}{n}
\bigg|\frac{\frac{\partial^{2}f_{1}}{\partial z_{j}\partial z_{l}}}{\frac{\partial f_{1}}{\partial z_{1}}}\bigg|
+\he{l=2}{n}\bigg|\frac{\frac{\partial f_{1}}{\partial z_{2}}}{\frac{\partial f_{2}}{\partial z_{2}}}\bigg|
\bigg|\frac{\frac{\partial^{2}f_{2}}{\partial z_{j}\partial z_{l}}}{\frac{\partial f_{1}}{\partial z_{1}}}\bigg|+\bigg|\frac{\frac{\partial f_{1}}{\partial z_{2}}}{\frac{\partial f_{1}}{\partial z_{1}}}\bigg|
\bigg|\frac{\frac{\partial f_{2}}{\partial z_{j}}}{\frac{\partial f_{2}}{\partial z_{2}}}\bigg|\bigg|
\frac{f''_{j}(z_{j})}{f'_{j}(z_{j})}\bigg|+\bigg|\frac{\frac{\partial f_{1}}{\partial z_{j}}}{\frac{\partial f_{1}}{\partial z_{1}}}\bigg|\bigg|
\frac{f''_{j}(z_{j})}{f'_{j}(z_{j})}\bigg|\bigg)\\
&&-p_{2}\bigg(\he{l=2}{n}\bigg|\frac{\frac{\partial^{2}f_{2}}{\partial z_{j}\partial z_{l}}}{\frac{\partial f_{2}}{\partial z_{2}}}\bigg|+\bigg|\frac{\frac{\partial f_{2}}{\partial z_{j}}}{\frac{\partial f_{2}}{\partial z_{2}}}\bigg|\bigg|
\frac{f''_{j}(z_{j})}{f'_{j}(z_{j})}\bigg|\bigg)\bigg]\geq 0.
\end{eqnarray*}
Thus it follows from Theorem A that $f\in K(D^{n}_{p})$.\hfill $\Box$

{\bf Remark 1}\quad Setting $f_{2}(z_{2},\cdots,z_{n})=f_{2}(z_{2})$ in Theorem 1, we get Theorem 2 of \cite{lz2009}.

Setting $n=3$ in Theorem 1, we get the following corollary, which gives an answer to Problem I for the case $n=3$.

{\bf Corollary 1}\quad Suppose that $p_j\ge 2,j=1,2,3$. Let $f_3:U\rightarrow C$ is holomorphic with $f_3(0)=0,f'_3(0)=1$, $f_{1}(z_{1},z_2,z_3), f_{2}(z_2,z_3):D_{p}^3\rightarrow C$ are holomorphic with $f_{1}(0,0,0)=0,\frac{\partial f_{1}}{\partial z_{1}}(0,0,0)=1,\frac{\partial f_{1}}{\partial z_{l}}(0,0,0)=0(l=2,3)$ and $f_{2}(0,0)=0,\frac{\partial f_{2}}{\partial z_{2}}(0,0)=1,\frac{\partial f_{2}}{\partial z_{3}}(0,0)=0$. Set
$$
f(z)=(f_{1}(z_{1},z_{2},z_3),f_{2}(z_{2},z_{3}),
f_{3}(z_{3})).
$$
If $f(z)$ satisfies the following conditions
\begin{eqnarray*}
&(1)&\frac{\partial f_{1}}{\partial z_{1}}\frac{\partial f_{2}}{\partial z_{2}}f'_3(z_3)\neq 0,\quad |z_3f''_3(z_3)|\leq|f'_3(z_3)|;\\
&(2)&\he{l=1}{3}\Big|z_{1}\frac{\partial^{2}f_{1}}{\partial z_{1}\partial z_{l}}\Big|\leq\Big|\frac{\partial f_{1}}{\partial z_{1}}\Big|,\quad \he{l=2}{3}\Big|z_{2}\frac{\partial^{2}f_{2}}{\partial z_{2}\partial z_{l}}\Big|\leq\Big|\frac{\partial f_{2}}{\partial z_{2}}\Big|;\\
&(3)&p_{1}\bigg(\he{l=1}{3}\bigg|\frac{\frac{\partial^{2}f_{1}}{\partial z_{2}\partial z_{l}}}{\frac{\partial f_{1}}{\partial z_{1}}}\bigg|+\he{l=2}{3}\bigg|\frac{\frac{\partial f_{1}}{\partial z_{2}}}{\frac{\partial f_{2}}{\partial z_{2}}}\bigg|\bigg|\frac{\frac{\partial^{2}f_{2}}{\partial z_{2}\partial z_{l}}}{\frac{\partial f_{1}}{\partial z_{1}}}\bigg|\bigg)\leq\bigg(1-\he{l=2}{3}\bigg|\frac{\frac{z_{2}\partial^{2}f_{2}}{\partial z_{2}\partial z_{l}}}{\frac{\partial f_{2}}{\partial z_{2}}}\bigg|\bigg)p_{2}|z_{2}|^{p_{2}-2}; \\
&(4)&p_{1}\bigg(\he{l=1}{3}\bigg|\frac{\frac{\partial^{2}f_{1}}{\partial z_{3}\partial z_{l}}}{\frac{\partial f_{1}}{\partial z_{1}}}\bigg|+\he{l=2}{3}\bigg|\frac{\frac{\partial f_{1}}{\partial z_{2}}}{\frac{\partial f_{2}}{\partial z_{2}}}\bigg|\bigg|\frac{\frac{\partial^{2}f_{2}}{\partial z_{3}\partial z_{l}}}{\frac{\partial f_{1}}{\partial z_{1}}}\bigg|+\bigg|\frac{\frac{\partial f_{1}}{\partial z_{2}}}{\frac{\partial f_{1}}{\partial z_{1}}}\bigg|\bigg|\frac{\frac{\partial f_{2}}{\partial z_{3}}}{\frac{\partial f_{2}}{\partial z_{2}}}\bigg|\bigg|\frac{f''_{3}(z_{3})}{f'_{3}(z_{3})}\bigg|+\bigg|\frac{\frac{\partial f_{1}}{\partial z_{3}}}{\frac{\partial f_{1}}{\partial z_{1}}}\bigg|\bigg|\frac{f''_{3}(z_{3})}{f'_{j}(z_{3})}\bigg|\bigg)\\
&&+p_{2}\bigg(\he{l=2}{n}\bigg|\frac{\frac{\partial^{2}f_{2}}{\partial z_{3}\partial z_{l}}}{\frac{\partial f_{2}}{\partial z_{2}}}\bigg|+\bigg|\frac{\frac{\partial f_{2}}{\partial z_{3}}}{\frac{\partial f_{2}}{\partial z_{2}}}\bigg|\bigg|\frac{f''_{3}(z_{3})}{f'_{3}(z_{3})}\bigg|\bigg)\leq p_{3}|z_{3}|^{p_{3}-2}\bigg(1-\bigg|z_{3}\frac{f''_{3}(z_{3})}{f'_{3}(z_{3})}\bigg|\bigg)
\end{eqnarray*}
for all $z=(z_1,z_2, z_3)\in D_p^3\setminus\{0\}$, then $f\in K(D_p^3)$.

Now let us give two examples to illustrate the application of Theorem 1 in the following.

{\bf Example 1.}\quad Suppose that $p_{j}\geq p_{1}\geq 2 (j=2,\cdots,n), 0< |\lambda|\leq 1$, and $k$ is a positive integer such that $k< \max\{p_{j}:j=1, \cdots, n\}\leq k+1$, let
$$
f(z)=(z_{1}+a_{1}z^{2}_{1}+\he{j=2}{n}a_{j}z^{k+1}_{j},z_{2}+a_{2}z^{2}_{2}+\he{j=3}{n}a_{j}z^{k+1}_{j},\frac{e^{\lambda z_{3}}-1}{\lambda},\cdots,\frac{e^{\lambda z_{n}}-1}{\lambda}),
$$
where $a=\max\{|a_{j}|: j=1,2,\cdots,n\}$. If $a\leq \frac{1-|\lambda|}{(k+1)^{2}+4}< \frac{1}{4}$ and
$$
\Big[\frac{p_{1}+p_2}{1-2a}+\frac{p_{1}(k+1)a}{(1-2a)^{2}}\Big]|a_{j}|\leq \frac{p_{j}(1-|\lambda|)}{(k+1)(k+|\lambda|)},\quad j=3,\cdots, n,
$$
then $f(z)\in K(D^{n}_{p})$.

{\bf Proof}\quad Let
$$
f_1(z)=z_{1}+a_{1}z^{2}_{1}+\he{j=2}{n}a_{j}z^{k+1}_{j},\, f_2(z)=z_{2}+a_{2}z^{2}_{2}+\he{j=3}{n}a_{j}z^{k+1}_{j}
$$
and $f_j(z)=\frac{e^{\lambda z_j}-1}{\lambda}$ for $j=3, \cdots, n$. Then it follows from $a=\max\{|a_{j}|: j=1, 2, \cdots, n\}<\frac{1}{4}$ that
\begin{eqnarray*}
|\frac{\partial f_{1}}{\partial z_{1}}|&=&|1+2a_{1}z_{1}|\geq 1-2a > 0,\quad|\frac{\partial f_{2}}{\partial z_{2}}|=|1+2a_{2}z_{2}|\geq 1-2a > 0,\\
\bigg|\frac{z_{1}\frac{\partial^{2} f_{1}}{\partial z_{1}\partial z_{l}}}{\frac{\partial f_{1}}{\partial z_{1}}}\bigg|&=&\frac{|2a_{1}z_{1}|}{|1+2a_{1}z_{1}|}\leq \frac{2a}{1-2a}\leq 1,\quad
\bigg|\frac{z_{2}\frac{\partial^{2} f_{2}}{\partial z_{2}\partial z_{l}}}{\frac{\partial f_{2}}{\partial z_{2}}}\bigg|=\frac{|2a_{2}z_{2}|}{|1+2a_{2}z_{2}|}\leq \frac{2a}{1-2a}\leq 1,\\
\bigg|\frac{z_{j}f''_{j}(z_{j})}{f'_{j}(z_{j})}\bigg|&=&|\lambda||z_{j}|\leq 1,\quad j=3,\cdots, n.
\end{eqnarray*}

By calculating straightforwardly, we also obtain
\begin{eqnarray*}
&& p_{2}|z_{2}|^{p_{2}-2}\bigg(1-\he{l=2}{n}\bigg|\frac{\frac{z_{2}\partial^{2} f_{2}}{\partial z_{2}\partial z_{l}}}{\frac{\partial f_{2}}{\partial z_{2}}}\bigg|\bigg)-p_{1}\bigg(\he{l=1}{n}\bigg|\frac{\frac{\partial^{2} f_{1}}{\partial z_{2}\partial z_{l}}}{\frac{\partial f_{1}}{\partial z_{1}}}\bigg|+\he{l=2}{n}\bigg|\frac{\frac{\partial f_{1}}{\partial z_{2}}}{\frac{\partial f_{2}}{\partial z_{2}}}\bigg|\bigg|\frac{\frac{\partial^{2} f_{2}}{\partial z_{2}\partial z_{l}}}{\frac{\partial f_{1}}{\partial z_{1}}}\bigg|\bigg)\\
&=&p_{2}|z_{2}|^{p_{2}-2}\bigg(1-\frac{|2a_{2}z_{2}|}{|1+2a_{2}z_{2}|}\bigg)-p_{1}\bigg(\frac{|k(k+1)a_{2}z^{k-1}_{2}|}
{|1+2a_{1}z_{1}|}+\frac{|(k+1)a_{2}z^{k}_{2}|}{|1+2a_{2}z_{2}|}\bigg|\frac{|2a_{2}|}{|1+2a_{1}z_{2}|}\bigg|\bigg)\\
&\geq &p_{2}|z_{2}|^{p_{2}-2}\bigg(1-\frac{2|a_{2}|}{1-2|a_{2}|}\bigg)-p_{1}\bigg(\frac{k(k+1)|a_{2}|}{1-2|a_{1}|}|z_{2}|^{p_{2}-2}
+\frac{(k+1)|a_{2}|}{1-2|a_{2}|}\frac{2|a_{2}|}{1-2|a_{1}|}|z_{2}|^{p_{2}-2}\bigg)\\
&\geq &p_{1}|z_{2}|^{p_{2}-2}\bigg(1-\frac{2a}{1-2a}-\frac{k(k+1)a}{1-2a}-\frac{(k+1)a}{1-2a}\bigg)\\
&=&\frac{p_{1}|z_{2}|^{p_{2}-2}}{1-2a}\{1-[(k+1)^{2}+4]a\}> 0.
\end{eqnarray*}

Since for $j=3,\cdots, n$, we have
\begin{eqnarray*}
\he{l=1}{n}\bigg|\frac{\frac{\partial^{2} f_{1}}{\partial z_{j}\partial z_{l}}}{\frac{\partial f_{1}}{\partial z_{1}}}\bigg|&=&\frac{|k(k+1)a_{j}z^{k-1}_{j}|}{|1+2a_{1}z_{1}|}\leq \frac{k(k+1)|a_{j}|}{1-2|a_{1}|}|z_{j}|^{k-1}\leq \frac{k(k+1)|a_{j}|}{1-2|a_{1}|}|z_{j}|^{p_{j}-2},\\
\he{l=2}{n}\bigg|\frac{\frac{\partial f_{1}}{\partial z_{2}}}{\frac{\partial f_{2}}{\partial z_{2}}}\bigg|\bigg|\frac{\frac{\partial^{2} f_{2}}{\partial z_{j}\partial z_{l}}}{\frac{\partial f_{1}}{\partial z_{1}}}\bigg|&=&\frac{|(k+1)a_{2}z^{k}_{2}|}{|1+2a_{2}z_{2}|}\frac{|k(k+1)a_j z_j^{k-1}|}{|1+2a_{1}z_{1}|}
\leq \frac{k(k+1)^{2}|a_{2}||a_{j}|}{(1-2|a_{2}|)(1-2|a_{1}|)}|z_j|^{p_{j}-2},\\
\bigg|\frac{\frac{\partial f_{1}}{\partial z_{2}}}{\frac{\partial f_{1}}{\partial z_{1}}}\bigg|
\bigg|\frac{\frac{\partial f_{2}}{\partial z_{j}}}{\frac{\partial f_{2}}{\partial z_{2}}}\bigg|\bigg|\frac{f''_{j}(z_{j})}{f'_{j}(z_{j})}\bigg|&=&\frac{|(k+1)a_{2}z^{k}_{2}|}{|1+2a_{1}z_{1}|}
\frac{|(k+1)a_{j}z^{k}_{j}|}{|1+2a_{2}z_{2}|}|\lambda|\leq \frac{(k+1)^{2}|\lambda||a_{2}||a_{j}|}{(1-2|a_{1}|)(1-2|a_{2}|)}|z_{j}|^{p_{j}-2},\\
\bigg|\frac{\frac{\partial f_1}{\partial z_{j}}}{\frac{\partial f_{1}}{\partial z_{1}}}\bigg|
\bigg|\frac{f''_{j}(z_{j})}{f'_{j}(z_{j})}\bigg|&=&\frac{|(k+1)a_{j}z^{k}_{j}|}{|1+2a_{1}z_{1}|}|\lambda|
\leq \frac{(k+1)|\lambda||a_{j}|}{1-2|a_{1}|}|z_{j}|^{p_{j}-2},\\
\he{l=2}{n}\bigg|\frac{\frac{\partial^{2} f_{2}}{\partial z_{j}\partial z_{l}}}{\frac{\partial f_{1}}{\partial z_{1}}}\bigg|&=&\frac{|k(k+1)a_{j}z^{k-1}_{j}|}{|1+2a_{1}z_{1}|}\leq \frac{k(k+1)|a_{j}|}{1-2|a_{1}|}|z_{j}|^{p_{j}-2},\\
\bigg|\frac{\frac{\partial f_{2}}{\partial z_{j}}}{\frac{\partial f_{2}}{\partial z_{2}}}\bigg|\bigg|
\frac{f''_{j}(z_{j})}{f'_{j}(z_{j})}\bigg|&=&\frac{|(k+1)a_{j}z^{k}_{j}|}{|1+2a_{2}z_{2}|}|\lambda|\leq
\frac{(k+1)|a_{j}||\lambda|}{1-2|a_{2}|}|z_{j}|^{p_{j}-2}.
\end{eqnarray*}

Therefore
\begin{eqnarray*}
&&p_{1}\bigg(\he{l=1}{n}\bigg|\frac{\frac{\partial^{2} f_{1}}{\partial z_{j}\partial z_{l}}}{\frac{\partial f_{1}}{\partial z_{1}}}\bigg|+\he{l=2}{n}\bigg|\frac{\frac{\partial f_{1}}{\partial z_{2}}}{\frac{\partial f_{2}}{\partial z_{2}}}\bigg|\bigg|\frac{\frac{\partial^{2} f_{2}}{\partial z_{j}\partial z_{l}}}{\frac{\partial f_{1}}{\partial z_{1}}}\bigg|+\bigg|\frac{\frac{\partial f_{1}}{\partial z_{2}}}{\frac{\partial f_{1}}{\partial z_{1}}}\bigg|\bigg|\frac{\frac{\partial f_{2}}{\partial z_{j}}}{\frac{\partial f_{2}}{\partial z_{2}}}\bigg|
\bigg|\frac{f''_{j}(z_{j})}{f'_{j}(z_{j})}\bigg|+\bigg|\frac{\frac{\partial f_{1}}{\partial z_{j}}}{\frac{\partial f_{1}}{\partial z_{1}}}\bigg|\bigg|\frac{f''_{j}(z_{j})}{f'_{j}(z_{j})}\bigg|\bigg)\\
&&+p_{2}\bigg(\he{l=2}{n}\bigg|
\frac{\frac{\partial^{2} f_{2}}{\partial z_{j}\partial z_{l}}}{\frac{\partial f_{1}}{\partial z_{1}}}\bigg|
+\bigg|\frac{\frac{\partial f_{2}}{\partial z_{j}}}{\frac{\partial f_{2}}{\partial z_{2}}}\bigg|\bigg|\frac{f''_{j}(z_{j})}{f'_{j}(z_{j})}\bigg|\bigg)\\
&\leq &p_{1}|z_{j}|^{p_{j}-2}\bigg(\frac{k(k+1)|a_{j}|}{1-2|a_{1}|}+\frac{k(k+1)^{2}|a_{2}||a_{j}|}
{(1-2|a_{1}|)(1-2|a_{2}|)}+\frac{(k+1)^{2}|\lambda||a_{2}||a_{j}|}{(1-2|a_{1}|)(1-2|a_{2}|)}\\
&&+\frac{(k+1)|a_{j}||\lambda|}{1-2|a_{1}|}\bigg)+p_{2}|z_{j}|^{p_{j}-2}\bigg(\frac{k(k+1)|a_{j}|}{1-2|a_{1}|}+\frac{(k+1)|\lambda||a_{j}|}{1-2|a_{2}|}\bigg)\\
&\leq & p_{1}|z_{j}|^{p_{j}-2}\bigg(\frac{k(k+1)|a_{j}|}{1-2a}+\frac{k(k+1)^{2}a|a_{j}|}{(1-2a)^{2}}+
\frac{(k+1)^{2}|\lambda|a|a_{j}|}{(1-2a)^{2}}+\frac{(k+1)|\lambda||a_{j}|}{1-2a}\bigg)\\
&& +p_{2}|z_{j}|^{p_{j}-2}\bigg(\frac{k(k+1)|a_{j}|}{1-2a}+\frac{(k+1)|\lambda||a_{j}|}{1-2a}\bigg)\\
&=&|z_{j}|^{p_{j}-2}(k+1)(k+|\lambda|)\bigg(\frac{p_{1}|a_{j}|}{1-2a}+\frac{p_{1}(k+1)a|a_{j}|}
{(1-2a)^{2}}+\frac{p_{2}|a_{j}|}{1-2a}\bigg)\\
&\leq &|z_{j}|^{p_{j}-2}(k+1)(k+|\lambda|)\frac{p_{j}(1-|\lambda|)}{(k+1)(k+\lambda)}=p_{j}|z_{j}|^{p_{j}-2}(1-|\lambda|)\\
&\leq &p_{j}|z_{j}|^{p_{j}-2}\bigg(1-\bigg|\frac{z_{j}f''_{j}(z_{j})}{f'_{j}(z_{j})}\bigg|\bigg).
\end{eqnarray*}
Hence it follows from Theorem 1 that $f\in K(D^{n}_{p})$.\hfill $\Box$

{\bf Example 2}\quad Suppose that $p_{j}\geq p_{1}\geq 2 (j=2,\cdots,n)$, and $k$ is a positive integer such that
$k< \max\{p_{j}:j=1, \cdots, n\}\leq k+1$, let
$$
f(z)=(z_{1}+a_{1}z^{2}_{1}+\he{j=2}{n}a_{j}z^{k+1}_{j},z_{2}+a_{2}z^{2}_{2}+\he{j=3}{n}a_{j}z^{2}_{3},z_{3}+a_{3}z^{2}_{3},\cdots,z_{n}+a_{n}z^{2}_{n})\\
$$
where $a=\max\{|a_{j}|: j=1,2,\cdots,n\}$. If $a\leq \frac{1}{(k+1)^{2}+4}< \frac{1}{4}$, and
$$
\bigg[\frac{p_{1}+p_{2}}{1-2a}+\frac{p_{1}(k+1)a}{(1-2a)^{2}}\bigg]|a_{j}|
\leq \frac{p_{j}(1-\frac{2|a_{j}|}{1-2|a_{j}|})}{(k+1)(k+\frac{2|a_{j}|}{1-2|a_{j}|})},\quad j=3,\cdots, n.
$$
Then $f(z)\in K(D^{n}_{p})$.

Next, we establish a sufficient condition for the biholomorphic convex mapping on $D_{p}^n$, which extend the main result of \cite{ll2014}.

{\bf Theorem 2}\quad Suppose that $n\geq 2,  p_j\ge 2,j=1,2,\cdots,n$. Let
$$
f(z)=(p_1(z_1,z_2,\cdots, z_n),  p_2(z_{2},z_{n}),\cdots ,
p_{n-1}(z_{n-1},z_{n}), p_n(z_{n})),
$$
where $z=(z_1, z_2,\cdots,z_n)\in D^n_p, p_n(z_n)\in
N(U), p_j(z_{j},z_{n}):D_p^2\to C$ is holomorphic with $p_j(0,0)=0, \frac{\partial p_j}{\partial
z_j}(0,0)=1, \frac{\partial p_j}{\partial
z_n}(0,0)=0$, and $p_1(z_1, \cdots, z_n):D^n_{p}\to C$ is holomorphic  with
$p_1(0,0,\cdots,0)=0, \frac{\partial p_1}{\partial
z_1}(0,0,\cdots,0)=1, \frac{\partial p_1}{\partial
z_l}(0,0,\cdots,0)=0$ for $2\leq l\leq n$. If $f(z)$ satisfies the following conditions
\begin{eqnarray*}
&(1)&\prod_{j=1}^{n-1}\frac{\partial p_{j}}{\partial z_{j}}p'_{n}(z_{n})\neq 0, \quad |z_{n}p''_{n}(z_{n})|\leq|p'_{n}(z_{n})|;\\
&(2)&\he{l=1}{n}\bigg|z_{1}\frac{\partial ^{2}p_{1}}{\partial z_{1}\partial z_{l}}\bigg|
\leq\bigg|\frac{\partial p_{1}}{\partial z_{1}}\bigg|,\quad
\bigg|z_{j}\frac{\partial^{2} p_{j}}{\partial z_{j}^{2}}\bigg|+\bigg|z_{j}\frac{\partial^{2} p_{j}}{\partial z_{j}\partial z_{n}}\bigg|\leq\bigg|\frac{\partial p_{j}}{\partial z_{j}}\bigg|(j=2,3,\cdots,n-1);\\
&(3)&\frac{p_1}{|\frac{\partial p_{1}}{\partial z_{1}}|}\bigg(\frac{|\frac{\partial p_{1}}{\partial z_{j}}|(|\frac{\partial^{2}p_{j}}{\partial z_{j}^{2}}|+|\frac{\partial^{2}p_{j}}{\partial z_{j}\partial z_{n}}|)}{|\frac{\partial p_{j}}{\partial z_{j}}|}+\he{l=1}{n}\bigg|\frac{\partial^{2}p_{1}}{\partial z_{j}\partial z_{l}}\bigg|\bigg)
\leq p_j|z_{j}|^{p_j-2}\bigg(1-\frac{\Big|z_{j}\frac{\partial^{2}p_{j}}{\partial z_{j}^{2}}\Big|+
\Big|z_{j}\frac{\partial^{2}p_{j}}{\partial z_{j}\partial z_{n}}\Big|}{\Big|\frac{\partial p_{j}}{\partial z_{j}}\Big|}\bigg)\\
&&\quad \quad (j=2,3,\cdots,n-1);\\
&(4)&\he{j=2}{n-1}\frac{p_{j}}{|\frac{\partial p_{j}}{\partial z_{j}}|}\bigg(\bigg|\frac{\partial^{2}p_{j}}{\partial z_{j}\partial z_{n}}\bigg|+\bigg|\frac{\partial^{2}p_{j}}{\partial z^{2}_{n}}\bigg|+\bigg|\frac{\partial p_{j}}{\partial z_{n}}\bigg|\bigg|\frac{p''_{n}(z_{n})}{p'_{n}(z_{n})}\bigg|\bigg)
+\frac{p_1}{\Big|\frac{\partial p_{1}}{\partial z_{1}}\Big|}\bigg(\he{l=1}{n}\bigg|\frac{\partial^{2}p_{1}}{\partial z_{n}\partial z_{l}}\bigg|\\
&&+\he{j=2}{n-1}\frac{\Big|\frac{\partial p_{1}}{\partial z_{j}}\Big|(\Big|\frac{\partial^{2}p_{j}}{\partial z_{j}\partial z_{n}}\Big|+\Big|\frac{\partial^{2} p_{j}}{\partial z_{n}^{2}}\Big|)}{\Big|\frac{\partial p_{j}}{\partial z_{j}}\Big|}
+\bigg|\frac{\partial p_{1}}{\partial z_{_{n}}}\bigg|\bigg|\frac{p''_{n}(z_{n})}{p'_{n}(z_{n})}\bigg|
+\he{j=2}{n-1}\Big|\frac{\partial p_{1}}{\partial z_{j}}\Big|
\bigg|\frac{\frac{\partial p_{j}}{\partial z_{n}}}{\frac{\partial p_{j}}{\partial z_{j}}}\bigg|\bigg|\frac{p''_{n}(z_{n})}{p'_{n}(z_{n})}\bigg|\bigg)\\
&&\leq p_n|z_{n}|^{p_n-2}\bigg(1-\bigg|\frac{z_{n}p''_{n}(z_{n})}{p'_{n}(z_{n})}\bigg|\bigg),
\end{eqnarray*}
for all $z=(z_1,\cdots, z_n)\in D_p^n\setminus\{0\}$, then $f\in K(D^{n}_{p})$.

{\bf Proof}\quad By calculating the Fr$\acute{e}$chet derivatives of $f(z)$ straightforwardly, we obtain
$$
Df(z)=\left (
\begin{array}{ccccc}
\frac{\partial p_1}{\partial z_1}& \frac{\partial p_1}{\partial
z_2}& \cdots &
 \frac{\partial p_1}{\partial z_{n-1}} &\frac{\partial p_1}{\partial z_n}\\
0 &\frac{\partial p_{2}}{\partial z_{2}}& \cdots &0 &\frac{\partial p_{2}}{\partial z_{n}}\\
\cdots &\cdots &\cdots &\cdots &\cdots \\
0 &0 &\cdots &\frac{\partial p_{n-1}}{\partial z_{n-1}}&\frac{\partial p_{n-1}}{\partial z_{n}}\\
0&0&\cdots &0&p_n'(z_n)
\end{array} \right ),
$$
$$
Df(z)^{-1}=\left (
\begin{array}{ccccc}
\frac{1}{\frac{\partial p_1}{\partial z_1}}& -\frac{\frac{\partial
p_1}{\partial z_2}} {\frac{\partial p_1}{\partial z_1}\frac{\partial p_{2}}{\partial z_{2}}}&
\cdots &-\frac{\frac{\partial p_1}{\partial z_{n-1}}}
{\frac{\partial p_1}{\partial z_1}\frac{\partial p_{n-1}}{\partial z_{n-1}}}
&-\frac{\frac{\partial p_1}{\partial z_n}} {\frac{\partial
p_1}{\partial z_1}p_n'(z_n)}+
\sum\limits_{j=2}^{n-1} \frac{\frac{\partial p_{j}}{\partial z_{n}}}{\frac{\partial p_{j}}{\partial z_{j}}}\frac{\frac{\partial p_1}{\partial z_j}}{\frac{\partial p_1}{\partial z_1}p_n'(z_n)}\\
0 & \frac{1}{\frac{\partial p_2}{\partial z_2}}& \cdots &0 &-\frac{\frac{\partial
p_2}{\partial z_n}} {\frac{\partial p_2}{\partial z_2}p_n'(z_n)}\\
\cdots &\cdots &\cdots &\cdots &\cdots \\
0 &0 &\cdots & \frac{1}{\frac{\partial p_{n-1}}{\partial z_{n-1}}}&-\frac{\frac{\partial
p_{n-1}}{\partial z_n}}{\frac{\partial p_{n-1}}{\partial z_{n-1}}p_n'(z_n)}\\
0&0&\cdots &0&\frac{1}{p_n'(z_n)}
\end{array} \right ),
$$
$$
\begin{array}{lll}
D^{2}f(z)(b,b)&=\left (
\begin{array}{ccccc}
\sum\limits_{l=1}^n\frac{\partial ^2p_1}{\partial z_1\partial
z_l}b_l& \sum\limits_{l=1}^n\frac{\partial ^2p_1}{\partial
z_2\partial z_l}b_l& \cdots & \sum\limits_{l=1}^n\frac{\partial
^2p_1}{\partial z_{n-1}\partial z_l}b_l &
C_{1}\\
0 &D_{2}& \cdots &0 & C_{2}\\
\cdots &\cdots &\cdots &\cdots &\cdots \\
0 &0 &\cdots &D_{n-1}&C_{n-1}\\
0&0&\cdots &0&C_{n}
\end{array} \right )\jz{b_{1}}{b_{2}}{\cdots }{b_{n-1}}{b_{n}}\\
&=\jz{\sum\limits_{j=1}^n\sum\limits_{l=1}^n\frac{\partial
^2p_1}{\partial z_j\partial z_l}b_lb_j}
{\frac{\partial^{2}p_{2}}{\partial z_{2}^{2}}b_{2}^{2}+2\frac{\partial^{2}p_{2}}{\partial z_{2}\partial z_{n}}b_{2}b_{n}+\frac{\partial^{2}p_{2}}{\partial z_{n}^{2}}b_{n}^{2}}{\cdots }
{\frac{\partial^{2}p_{n-1}}{\partial z_{n-1}^{2}}b_{n-1}^{2}+2\frac{\partial^{2}p_{n-1}}{\partial z_{n-1}\partial z_{n}}b_{n-1}b_{n}+\frac{\partial^{2}p_{n-1}}{\partial z_{n}^{2}}b_{n}^{2}}{p_n''(z_n)b_n^2},
\end{array}
$$
where
\begin{eqnarray*}
C_{1}&=&\sum\limits_{l=1}^n\frac{\partial ^2p_1}{\partial z_n\partial z_l}b_l,\\
C_{j}&=&\frac{\partial^{2}p_{j}}{\partial z_{n}\partial z_{j}}b_j+ \frac{\partial^{2}p_{j}}{\partial z_{n}^{2}}b_n\, (j=2,3,\cdots,n-1),\quad C_{n}=p_n''(z_n)b_n,\\
D_{j}&=&\frac{\partial^{2}p_{j}}{\partial z_{j}^{2}}b_j+ \frac{\partial^{2}p_{j}}{\partial z_{j}\partial z_{n}}b_n\, (j=2,3,\cdots,n-1).
\end{eqnarray*}

Taking $z=(z_{1}, \cdots ,z_{n})\in D^n_p\setminus\{0\}, b=(b_{1},\cdots, b_{n})\in C^n$ such that
$$
\mbox{Re}\bigg\{\he{j=1}{n}p_{j}\bigg|\frac{z_{j}}{\rho(z)}\bigg|^{p_{j}}\frac{b_{j}}{z_{j}}\bigg\}=0,
$$
by the hypothesis of Theorem 2, we have
\begin{eqnarray*}
&&J_{f}(z,b)\ge \he{j=1}{n}\frac{p_{j}|z_{j}|^{p_{j}-2}}{\rho(z)^{p_{j}}}|b_{j}|^{2}-
2\he{j=1}{n}\frac{p_{j}}{\rho(z)}\Big|\frac{z_{j}}{\rho(z)}\Big|^{p_{j}}
\mbox{Re}\langle Df(z)^{-1}D^{2}f(z)(b,b), \frac{\partial \rho}{\partial\overline{z}}\rangle\\
&=&\he{j=1}{n}\frac{p_{j}|z_{j}|^{p_{j}-2}}{\rho(z)^{p_{j}}}|b_{j}|^{2}
-\mbox{Re}\bigg\{\frac{1}{\frac{\partial p_{1}}{\partial z_{1}}}\bigg[\he{j=1}{n}\he{l=1}{n}\frac{\partial ^2p_1}{\partial z_j\partial z_l}b_l b_{j}
-\he{j=2}{n-1}\frac{\frac{\partial p_{1}}{\partial z_{j}}}{\frac{\partial p_{j}}{\partial z_{j}}}\bigg(\frac{\partial^{2}p_{j}}{\partial z_{j}^{2}}b_{j}^{2}
+2\frac{\partial^{2}p_{j}}{\partial z_{j}\partial z_{n}}b_{j}b_{n}\\
&& +\frac{\partial^{2}p_{j}}{\partial z_{n}^{2}}b_{n}^{2}\bigg)
-\frac{\partial p_{1}}{\partial z_{n}}\frac{p''_{n}(z_{n})}{p'_{n}(z_{n})} b_{n}^{2}
+\he{j=2}{n-1}\frac{\frac{\partial p_{j}}{\partial z_{n}}\frac{\partial p_{1}}{\partial z_{j}}}{\frac{\partial p_{j}}{\partial z_{j}}\frac{\partial p_{1}}{\partial z_{1}}}\frac{p''_{n}(z_{n})}{p'_{n}(z_{n})}b_{n}^{2}\bigg]\frac{p_{1}|z_{1}|
^{p_{1}}}{z_{1}\rho(z)^{p_{1}}}\\
&& +\he{j=2}{n-1}\frac{1}{\frac{\partial p_{j}}{\partial z_{j}}}\bigg(\frac{\partial^{2}p_{j}}{\partial z_{j}^{2}}b_{j}^{2}+2\frac{\partial^{2}p_{j}}{\partial z_{j}\partial z_{n}}b_{j}b_{n}+\frac{\partial^{2}p_{j}}{\partial z_{n}^{2}}b_{n}^{2}-\frac{\partial p_{j}}{\partial z_{ n}} \frac{p''_{n}(z_{n})}{p'_{n}(z_{n})}b^{2}_{n}\bigg)\frac{p_{j}|z_{j}|
^{p_{j}}}{z_{j}\rho(z)^{p_{j}}}
+\frac{p''_{n}(z_{n})}{p'_{n}(z_{n})}b_{n}^{2}\frac{p_{n}|z_{n}|^{p_{n}}}
{z_{n}\rho(z)^{p_{n}}}\bigg\}\\
&\geq &\he{j=1}{n}\frac{p_{j}|z_{j}|^{p_{j}-2}}{\rho(z)^{p_{j}}}|b_{j}|^{2}-
\frac{1}{|\frac{\partial p_{1}}{\partial z_{1}}|}\bigg[\he{j=1}{n}\he{l=1}{n}\bigg|\frac{\partial ^2p_1}{\partial z_j\partial z_l}\bigg||b_{j}|^{2}
+\he{j=2}{n-1}\bigg|\frac{\frac{\partial p_{1}}{\partial z_{j}}}{\frac{\partial p_{j}}{\partial z_{j}}}\bigg|\bigg(\bigg|\frac{\partial^{2}p_{j}}{\partial z_{j}^{2}}\bigg||b_{j}|^{2}
+\bigg|\frac{\partial^{2}p_{j}}{\partial z_{j}\partial z_{n}}\bigg||b_{j}|^{2}\\
&& +\bigg|\frac{\partial^{2}p_{j}}{\partial z_{j}\partial z_{n}}\bigg||b_{n}|^{2}
+\bigg|\frac{\partial^{2}p_{j}}{\partial z_{n}^{2}}\bigg||b_{n}|^{2}\bigg)
+\bigg|\frac{\partial p_{1}}{\partial z_{n}}\bigg|\bigg|\frac{p''_{n}(z_{n})}{p'_{n}(z_{n})}\bigg||b_{n}|^{2}
+\he{j=2}{n-1}\bigg|\frac{\frac{\partial p_{j}}{\partial z_{n}}\frac{\partial p_{1}}{\partial z_{j}}}{\frac{\partial p_{j}}{\partial z_{j}}} \bigg|\bigg|\frac{p''_{n}(z_{n})}{p'_{n}(z_{n})}\bigg||b_{n}|^{2}\bigg]
\frac{p_1|z_{1}|^{p_1-1}}{\rho(z)^{p_1}}\\
\end{eqnarray*}
\begin{eqnarray*}
&& -\he{j=2}{n-1}\frac{1}{|\frac{\partial p_{j}}{\partial z_{j}}|}\bigg(\bigg|\frac{\partial^{2}p_{j}}{\partial z_{j}^{2}}\bigg||b_{j}|^{2}+\bigg|\frac{\partial^{2}p_{j}}{\partial z_{j}\partial z_{n}}\bigg||b_{j}|^{2}
+\bigg|\frac{\partial^{2}p_{j}}{\partial z_{j}\partial z_{n}}\bigg||b_{n}|^{2}\\
&&+\bigg|\frac{\partial^{2}p_{j}}{\partial z_{n}^{2}}\bigg||b_{n}|^{2}+\bigg|\frac{\partial p_{j}}{\partial z_{ n}}\bigg| \bigg|\frac{p''_{n}(z_{n})}{p'_{n}(z_{n})}\bigg||b_{n}|^{2}\bigg)
\frac{p_{j}|z_{j}|^{p_{j}-1}}{\rho(z)^{p_{j}}}
-\bigg|\frac{p''_{n}(z_{n})}{p'_{n}(z_{n})}\bigg||b_{n}|^{2}\frac{p_n|z_{n}|^{p_n-1}}{\rho(z)^{p_n}}\\
&=&|b_{1}|^2\frac{p_1|z_{1}|^{p_1-2}}{\rho(z)^{p_1}}\bigg(1-\frac{\he{l=1}{n}|z_{1}\frac{\partial^{2}p_{1}}{\partial z_{1}\partial z_{l}}|}{|\frac{\partial p_{1}}{\partial z_{1}}|}\bigg)
+\he{j=2}{n-1}|b_{j}|^{2}\bigg[\frac{p_j|z_{j}|^{p_j-2}}{\rho(z)^{p_j}}\bigg(1-\frac{|{z_{j}\frac{\partial^{2}p_{j}}{\partial z_{j}^{2}}}|+|z_{j}\frac{\partial^{2}p_{j}}{\partial z_{j}\partial z_{n}}|}{|\frac{\partial p_{j}}{\partial z_{j}}|}\bigg)\\
&&-\frac{p_1|z_{1}|^{p_{1}-1}}{\rho(z)^{p_{1}}\bigg|\frac{\partial p_{1}}{\partial z_{1}}\bigg|}\bigg(\he{l=1}{n}|\frac{\partial^{2}p_{1}}{\partial z_{j}\partial z_{l}}|+\bigg|\frac{\partial p_{1}}{\partial z_{j}}\bigg| \frac{|\frac{\partial^{2}p_{j}}{\partial z_{j}^{2}}|
+|\frac{\partial^{2}p_{j}}{\partial z_{j}\partial z_{n}}|}{|\frac{\partial p_{j}}{\partial z_{j}}|}\bigg)\bigg]
+|b_{n}|^{2}\bigg[\frac{p_n|z_{n}|^{p_{n}-2}}{\rho(z)^{p_n}}(1-|\frac{z_{n}p''_{n}(z_{n})}{p'_{n}(z_{n})}|)\\
&&-\frac{p_1|z_{1}|^{p_{1}-1}}{\rho(z)^{p_{1}}|\frac{\partial p_{1}}{\partial z_{1}}|}\bigg(\he{l=1}{n}\bigg|\frac{\partial^{2}p_{1}}{\partial z_{n}\partial z_{l}}\bigg|+\sum\limits_{j=2}^{n-1}\bigg|\frac{\partial p_{1}}{\partial z_{j}}\bigg|\frac{|\frac{\partial^{2}p_{j}}{\partial z_{n}^{2}}|
+|\frac{\partial^{2}p_{j}}{\partial z_{j}\partial z_{n}}|}{|\frac{\partial p_{j}}{\partial z_{j}}|} +\bigg|\frac{\partial p_{1}}{\partial z_{n}}\bigg|\bigg|\frac{p''_{n}(z_{n})}{p'_{n}(z_{n})}\bigg|
+\he{j=2}{n-1}\bigg|\frac{\partial p_{1}}{\partial z_{j}}\bigg|\bigg| \frac{\frac{\partial p_{j}}{\partial z_{n}}}{\frac{\partial p_{j}}{\partial z_{j}}}\bigg|\bigg|\frac{p''_{n}(z_{n})}{p'_{n}(z_{n})}\bigg|\bigg)\\
&& -\he{j=2}{n-1}\bigg|\frac{\frac{\partial^{2}p_{j}}{\partial z_{j}\partial z_{n}}}{\frac{\partial p_{j}}{\partial z_{j}}}\bigg|\frac{p_j|z_{j}|^{p_j-1}}{\rho(z)^{p_j}}
-\he{j=2}{n-1}\bigg|\frac{\frac{\partial^{2}p_{j}}{\partial z_{n}^{2}}}{\frac{\partial p_{j}}{\partial z_{j}}}\bigg|\frac{p_j|z_{j}|^{p_j-1}}{\rho(z)^{p_j}}
-\he{j=2}{n-1}\bigg|\frac{\frac{\partial p_{j}}{\partial z_{n}}}{\frac{\partial p_{j}}{\partial z_{j}}}\bigg|\bigg|\frac{p''_{n}(z_{n})}{p'_{n}(z_{n})}\bigg|\frac{p_j|z_{j}|^{p_j-1}}{\rho(z)^{p_j}}\bigg]\\
&\geq &|b_{1}|^2\frac{p_1|z_{1}|^{p_1-2}}{\rho(z)^{p_1}}\bigg(1-\frac{\he{l=1}{n}|z_{1}\frac{\partial^{2}p_{1}}{\partial z_{1}\partial z_{l}}|}{|\frac{\partial p_{1}}{\partial z_{1}}|}\bigg)
+\he{j=2}{n-1}|b_{j}|^{2}\bigg[p_j|z_{j}|^{p_j-2}\bigg(1-\frac{|{z_{j}\frac{\partial^{2}p_{j}}{\partial z_{j}^{2}}}|+|z_{j}\frac{\partial^{2}p_{j}}{\partial z_{j}\partial z_{n}}|}{|\frac{\partial p_{j}}{\partial z_{j}}|}\bigg)\\
&& -\frac{p_1}{|\frac{\partial p_{1}}{\partial z_{1}}|}\bigg(\he{l=1}{n}|\frac{\partial^{2}p_{1}}{\partial z_{j}\partial z_{l}}|+\bigg|\frac{\partial p_{1}}{\partial z_{j}}\bigg| \frac{|\frac{\partial^{2}p_{j}}{\partial z_{j}^{2}}|
+|\frac{\partial^{2}p_{j}}{\partial z_{j}\partial z_{n}}|}{|\frac{\partial p_{j}}{\partial z_{j}}|}\bigg)\bigg]
+|b_{n}|^{2}\bigg[p_n|z_{n}|^{p-2}(1-|\frac{z_{n}p''_{n}(z_{n})}{p'_{n}(z_{n})}|)\\
&& -\frac{p_1}{|\frac{\partial p_{1}}{\partial z_{1}}|}\bigg(\he{l=1}{n}\bigg|\frac{\partial^{2}p_{1}}{\partial z_{n}\partial z_{l}}\bigg|+\sum\limits_{j=2}^{n-1}\bigg|\frac{\partial p_{1}}{\partial z_{j}}\bigg|\frac{|\frac{\partial^{2}p_{j}}{\partial z_{n}^{2}}|
+|\frac{\partial^{2}p_{j}}{\partial z_{j}\partial z_{n}}|}{|\frac{\partial p_{j}}{\partial z_{j}}|} +\bigg|\frac{\partial p_{1}}{\partial z_{n}}\bigg|\bigg|\frac{p''_{n}(z_{n})}{p'_{n}(z_{n})}\bigg|
+\he{j=2}{n-1}\bigg|\frac{\partial p_{1}}{\partial z_{j}}\bigg|\bigg| \frac{\frac{\partial p_{j}}{\partial z_{n}}}{\frac{\partial p_{j}}{\partial z_{j}}}\bigg|\bigg|\frac{p''_{n}(z_{n})}{p'_{n}(z_{n})}\bigg|\bigg)\\
&& -\he{j=2}{n-1}\bigg|\frac{\frac{\partial^{2}p_{j}}{\partial z_{j}\partial z_{n}}}{\frac{\partial p_{j}}{\partial z_{j}}}\bigg|p_j
-\he{j=2}{n-1}\bigg|\frac{\frac{\partial^{2}p_{j}}{\partial z_{n}^{2}}}{\frac{\partial p_{j}}{\partial z_{j}}}\bigg|p_j
-\he{j=2}{n-1}\bigg|\frac{\frac{\partial p_{j}}{\partial z_{n}}}{\frac{\partial p_{j}}{\partial z_{j}}}\bigg|\bigg|\frac{p''_{n}(z_{n})}{p'_{n}(z_{n})}\bigg|p_j\bigg]\\
&\geq & 0.
\end{eqnarray*}
Thus it follows from Theorem A that $f\in K(D^{n}_{p})$. The proof is complete.\hfill $\Box$

{\bf Remark 2}\quad Setting $p_{j}(z_{j},z_{n})=f_{j}(z_{j})+p_j(z_n)\, j=2,\cdots, n$ in Theorem 2, we get Theorem 2.1 of \cite{ll2014}. Setting $n=3$ in Theorem 2, we get Corollary 1.

{\bf Example 3}\quad Suppose that $p_{j}\geq p_{1}\geq 2 (j=2,\cdots,n),0<|\lambda|\leq1$, and $k$ is a positive integer such that $k<\max\{p_{j}:j=1, \cdots, n\}\leq k+1$. Let
$$p(z)=(z_{1}+\he{j=2}{n-1}a_{j}z^{k+1}_{j}+a_{n}z_{1}z^{k+1}_{n}, z_{2}+a_{2}z_{2}z^{k+1}_{n},\cdots,z_{n-1}+a_{n-1}
z_{n-1}z^{k+1}_{n},\frac{e^{\lambda z_{n}-1}}{\lambda}),$$
where $a=\max\bigg\{|a_{j}|:j=2,\cdots,n\bigg\}$. If $a\leq\frac{1-|\lambda|}{2(k+1)^{2}(k+1+|\lambda|)+1+|\lambda|}<1$, and
$$
\he{j=2}{n-1}\frac{p_{j}|a_{j}|}{1-|a_{j}|}+\frac{a\, p_1}{1-a}\bigg(1+(k+1)\he{j=2}{n-1}\frac{|a_{j}|}{1-|a_{j}|}\bigg)
\leq\frac{p_{n}(1-|\lambda|)}{(k+1)(k+1+|\lambda|)},
$$
then $p(z)\in K(D^{n}_{p})$.

{\bf Proof}\quad Put
\begin{eqnarray*}
p_{1}(z_{1},z_{2},\cdots,z_{n})&=&z_{1}+\he{j=2}{n-1}a_{j}z^{k+1}_{j}+a_{n}z_{1}z^{k+1}_{n},\\
p_{j}(z_{j},z_{n})&=&z_{j}+a_{j}z_{j}z^{k+1}_{n}\, (j=2,3,\cdots,n-1), \quad p_{n}(z_{n})
=\frac{e^{\lambda z_{n}}-1}{\lambda}.
\end{eqnarray*}
Then it follows from $a=\max\bigg\{|a_{j}|:j=2,\cdots,n\bigg\}\leq\frac{1-|\lambda|}{2(k+1)^{2}(k+1+|\lambda|)+1-|\lambda|}<\frac{1}{k+2}<1$ that
\begin{eqnarray*}
\bigg|\frac{\partial p_{1}}{\partial z_{1}}\bigg|&=&|1+a_{n}z^{k+1}_{n}|\geq 1-|a_{n}|\geq 1-a>0,\\
\bigg|\frac{\partial p_{j}}{\partial z_{j}}\bigg|&=&|1+a_{j}z^{k+1}_{n}|\geq 1-|a_{j}|\geq 1-a>0,\\
p'_{n}(z_{n})&=&e^{\lambda z_{n}}\neq 0,\\
\bigg|\frac{\partial p_{1}}{\partial z_{1}}\bigg|-\he{l=1}{n}\bigg|z_{1}\frac{\partial^{2}p_{1}}{\partial z_{1}\partial z_{l}}\bigg|&=&|1+a_{n}z^{k+1}_{n}|-|(k+1)a_{n}z_{1}z^{k}_{n}|\geq 1-(k+2)|a_{n}|>0,\\
\bigg|\frac{\partial p_{j}}{\partial z_{j}}\bigg|-\bigg|z_{j}\frac{\partial^{2}p_{j}}{\partial z^{2}_{j}}\bigg|
-\bigg|\frac{z_{j}\partial^{2}p_{j}}{\partial z_{j}\partial z_{n}}\bigg|&=&|1+a_{j}z^{k+1}_{n}|-(k+1)|a_{j}z_{j}z^{k}_{n}|\geq 1-(k+2)|a_{j}|>0,
\end{eqnarray*}
and $\bigg|\frac{z_{n}p''_{n}(z_{n})}{p'_{n}(z_{n})}\bigg|=|\lambda||z|\leq |\lambda|< 1$.
By calculating straightforwardly, we also obtain
\begin{eqnarray*}
&&p_{j}|z_{j}|^{p_{j}-2}
\bigg(1-\frac{\bigg|z_{j}\frac{\partial^{2}p_{j}}{\partial z^{2}_{j}}\bigg|+\bigg|z_{j}\frac{\partial^{2}p_{j}}{\partial z_{j}\partial z_{n}}\bigg|}{\bigg|\frac{\partial p_{j}}{\partial z_{j}}\bigg|}\bigg)-\frac{p_{1}}{|\frac{\partial p_{1}}{\partial z_{1}}|}
\bigg(\frac{\bigg|\frac{\partial p_{1}}{\partial z_{j}}\bigg|\bigg(\bigg|\frac
{\partial^{2}p_{j}}{\partial z^{2}_{j}}\bigg|+\bigg|\frac{\partial^{2}p_{j}}{\partial z_{j}\partial z_{n}}\bigg|\bigg)}{\bigg|\frac{\partial p_{j}}{\partial z_{j}}\bigg|}+\he{l=1}{n}\bigg|\frac{\partial^{2}p_{1}}{\partial z_{j}\partial z_{l}}\bigg|\bigg)\\
&=&p_{j}|z_{j}|^{p_{j}-2}\bigg(1-\frac{|(k+1)z_{j}a_{j}z^{k}_{n}|}{|1+a_{j}z^{k+1}_{n}|}\bigg)
-\frac{p_{1}}{|1+a_{n}z^{k+1}_{n}|}\bigg(\frac{|(k+1)a_{j}z^{k}_{j}||(k+1)a_{j}z^{k}_{n}|}{|1+a_{j}z^{k+1}_{n}|}
+|k(k+1)a_{j}z^{k-1}_{j}|\bigg)\\
&\geq & p_{j}|z_{j}|^{p_{j}-2}\bigg(1-\frac{(k+1)|a_{j}|}{1-|a_{j}|}\bigg)-\frac{p_{j}}{1-|a_{n}|}\bigg(\frac{(k+1)^{2}|a_{j}|^{2}}{1-|a_{j}|}+k(k+1)|a_{j}|\bigg)|z_{j}|^{k-1}\\
&\geq &p_{j}|z_{j}|^{p_{j}-2}\bigg(1-\frac{(k+1)a}{1-a}\bigg)-\frac{p_{j}}{1-a}\bigg(\frac{(k+1)^{2}a^{2}}{1-a}+(k+1)^{2}a\bigg)|z_{j}|^{p_{j}-2}\\
&=&\frac{p_{j}|z_{j}|^{p_{j}-2}}{1-a}\bigg(1-(k+2)a-\frac{(k+1)^{2}a}{1-a}\bigg)\\
&\geq&\frac{p_{j}|z_{j}|^{p_{j}-2}}{1-a}\bigg(1-\frac{(k+2)(1-|\lambda|)}{2(k+1)^{2}(k+1+|\lambda|)+1-|\lambda|}-\frac{1-|\lambda|}{2(k+1+|\lambda|)}\bigg)\\
\end{eqnarray*}
\begin{eqnarray*}
&\geq&\frac{p_{j}|z_{j}|^{p_{j}-2}}{1-a}\bigg[1-\frac{[(k+2)+(k+1)^{2}](1-|\lambda|)}{2(k+1)^{2}(k+1+|\lambda|)}\bigg]>0,\quad j=2, \cdots, n-1.
\end{eqnarray*}

Since
\begin{eqnarray*}
&&\he{j=2}{n-1}p_j\bigg|\frac{\frac{\partial^{2}p_{j}}{\partial z_{j}\partial z_{n}}}{\frac{\partial p_{j}}{\partial z_{j}}}\bigg|=\he{j=2}{n-1}p_j\bigg|\frac{(k+1)a_{j}z^{k}_{n}}{1+a_{j}z^{k+1}_{n}}\bigg|\leq(k+1)
\he{j=2}{n-1}\frac{p_j|a_{j}|}{1-|a_{j}|}|z_{n}|^{p_{n}-2},\\
&&\he{j=2}{n-1}p_j\bigg|\frac{\frac{\partial^{2}p_{j}}{\partial z^{2}_{n}}}{\frac{\partial p_{j}}{\partial z_{j}}}\bigg|=\he{j=2}{n-1}p_j\bigg|\frac{k(k+1)a_{j}z_{j}z^{k-1}_{n}}{1+a_{j}z^{k+1}_{n}}\bigg|\leq k(k+1)\he{j=2}{n-1}\frac{p_j|a_{j}|}{1-|a_{j}|}|z_{n}|^{p_{n}-2},\\
&&\he{j=2}{n-1}p_j\bigg|\frac{\frac{\partial p_{j}}{\partial z_{n}}}{\frac{\partial p_{j}}{\partial z_{j}}}\bigg|\bigg|\frac{p''_{n}(z_{n})}{p'_{n}(z_{n})}\bigg|=\he{j=2}{n-1}p_j\bigg|\frac{(k+1)a_{j}z_{j}z^{k}_{n}}{1+a_{j}z^{k+1}_{n}}\bigg|
|\lambda|\leq |\lambda|(k+1)\he{j=2}{n-1}\frac{p_j|a_{j}|}{1-|a_{j}|}|z_{n}|^{p_{n}-2},
\end{eqnarray*}
we have
\begin{eqnarray*}
&&\he{j=2}{n-1}\frac{p_{j}}{|\frac{\partial p_{j}}{\partial z_{j}}|}\bigg(\bigg|\frac{\partial^{2}p_{j}}{\partial z_{j}\partial z_{n}}\bigg|+\bigg|\frac{\partial^{2}p_{j}}{\partial z^{2}_{n}}\bigg|+\bigg|\frac{\partial p_{j}}{\partial z_{n}}\bigg|\bigg|\frac{p''_{n}(z_{n})}{p'_{n}(z_{n})}\bigg|\bigg)+\frac{p_{1}}{|\frac{\partial p_{1}}{\partial z_{1}}|}\bigg(\he{l=1}{n}\bigg|\frac{\partial^{2}p_{1}}{\partial z_{n}\partial z_{l}}\bigg|\\
&&+\he{j=2}{n-1}\frac{|\frac{\partial p_{1}}{\partial z_{j}}|\bigg(\bigg|\frac{\partial^{2}p_{j}}{\partial z_{j}\partial z_{n}}\bigg|+\bigg|\frac{\partial^{2}p_{j}}{\partial z^{2}_{n}}\bigg|\bigg)}{|\frac{\partial p_{j}}{\partial z_{j}}|}+\bigg|\frac{\partial p_{1}}{\partial z_{n}}\bigg|\bigg|\frac{p''_{n}(z_{n})}{p'_{n}(z_{n})}\bigg|+\he{j=2}{n-1}\bigg|\frac{\partial p_{1}}{\partial z_{j}}\bigg|\bigg|\frac{\frac{\partial p_{j}}{\partial z_{n}}}{\frac{\partial p_{j}}{\partial z_{j}}}\bigg|\bigg|\frac{p''_{n}(z_{n})}{p'_{n}(z_{n})}\bigg|\bigg)\\
&\leq & |z_{n}|^{p_{n}-2}\bigg[(k+1)\he{j=2}{n-1}\frac{p_j|a_{j}|}{1-|a_{j}|}+k(k+1)\he{j=2}{n-1}\frac{p_j|a_{j}|}{1-|a_{j}|}
+|\lambda|(k+1)\he{j=2}{n-1}\frac{p_j|a_{j}|}{1-|a_{j}|}\bigg]\\
&&+\frac{p_{1}}{1-|a_{n}|}\bigg[(k+1)|a_{n}|+k(k+1)|a_{n}|+\he{j=2}{n-1}\frac{(k+1)|a_{j}|[(k+1)|a_{j}|+k(k+1)|a_{j}|]}{1-|a_{j}|}\\
&&+(k+1)|a_{n}||\lambda|+\he{j=2}{n-1}
(k+1)|a_{j}|\frac{(k+1)|a_{j}||\lambda|}{1-|a_{j}|}\bigg]|z_{n}|^{k-1}\\
&\leq & |z_{n}|^{p_{n}-2}(k+1)(k+1+|\lambda|)\he{j=2}{n-1}\frac{p_j|a_{j}|}{1-|a_{j}|}+\frac{p_{1}}{1-a}\bigg[(k+1)a+k(k+1)a\\
&&+\he{j=2}{n-1}\frac{(k+1)|a_{j}|(k+1)^{2}|a_{j}|}{1-|a_{j}|}
+(k+1)a|\lambda|+\he{j=2}{n-1}\frac{(k+1)|a_{j}|(k+1)|a_{j}||\lambda|}{1-|a_{j}|}\bigg]|z_{n}|^{p_{n}-2}\\
&\leq & |z_{n}|^{p_{n}-2}(k+1)(k+1+|\lambda|)\he{j=2}{n-1}\frac{p_j|a_{j}|}{1-|a_{j}|}\\
&&+\frac{p_{1}}{1-a}(k+1)(k+1+|\lambda|)a|z_{n}|^{p_{n}-2}
\bigg[1+(k+1)\he{j=2}{n-1}\frac{|a_{j}|}{1-|a_{j}|}\bigg]\\
&\leq & |z_{n}|^{p_{n}-2}(k+1)(k+1+|\lambda|)\bigg[\he{j=2}{n-1}\frac{p_{j}|a_{j}|}{1-|a_{j}|}+\frac{a\, p_1}{1-a}\bigg(1+(k+1)\he{j=2}{n-1}\frac{|a_{j}|}{1-|a_{j}|}\bigg)\bigg]\\
&\leq& p_{n}|z_{n}|^{p_{n}-2}(1-|\lambda|)\leq p_{n}|z_{n}|^{p_{n}-2}\bigg(1-\bigg|\frac{z_{n}p''_{n}(z_{n})}{p'_{n}(z_{n})}\bigg|\bigg),
\end{eqnarray*}
hence by Theorem 2, we obtain that $f \in K(D^{n}_{p})$. \hfill $\Box$

By applying the same method of the proof for example 3, we only need to let $\frac{2|a_{n}|}{1-2|a_{n}|}$ instead of $|\lambda|$, we may prove the following result.

{\bf Example 4}\quad  Suppose that $p_{j}\geq p_{1}\geq 2 (j=2,\cdots,n),0<|a_{n}|\leq\frac{1}{4}$, and $k$ is a positive integer such that $k<\max\{p_{j}:j=1, \cdots, n\}\leq k+1$. Let
$$
p(z)=(z_{1}+\he{j=2}{n-1}a_{j}z^{k+1}_{j}+a_{n}z_{1}z^{k+1}_{n},z_{2}+a_{2}z_{n}z^{k+1}_{n},\cdots,z_{n-1}+a_{n-1}
z_{n-1}z^{k+1}_{n},z_{n}+a_{n}z^{2}_{n}),
$$
where $a=\max\bigg\{|a_{j}|:j=2,\cdots,n\bigg\}$. If $$
a\leq\frac{1-\frac{2|a_{n}|}{1-2|a_{n}|}}{2(k+1)^{2}(k+1+\frac{2|a_{n}|}{1-2|a_{n}|})+1-\frac{2|a_{n}|}{1-2|a_{n}|}}<1,
$$
and
$$\he{j=2}{n-1}\frac{p_{j}|a_{j}|}{1-|a_{j}|}+\frac{p_{1}a}{1-a}\bigg[1+(k+1)\he{j=2}{n-1}\frac{|a_{j}|}{1-|a_{j}|}\bigg]
\leq\frac{p_{n}(1-\frac{2|a_{n}|}{1-2|a_{n}|})}{(k+1)(k+1+\frac{2|a_{n}|}{1-2|a_{n}|})},$$
then $p(z)\in K(D^{n}_{p})$.

Applying the same method as our proof of Theorem 1, we may show the following theorem.

{\bf Theorem 3}\quad Suppose that $n\geq 2, p_j\ge 2, j=1,2,\cdots,n$, $k$ is a
positive integer such that $2\leq k\leq n$, and $f_{j},\, p_j:U\rightarrow C$ is holomorphic with $f_{j}(0)=0,f'_{j}(0)=0(j=2,3,\cdots,n),f_{k}(z_{k})=0$, $p_{1}(z_{1},z_{2},\cdots,z_{n}): D_{p}\rightarrow C$ is holomorphic with $p_{1}(0,0,\cdots,0)=0,\frac{\partial p_{1}}{\partial z_{1}}(0,0,\cdots,0)=1,\frac{\partial p_{1}}{\partial z_{l}}(0,0,\cdots,0)=0(l=2,3,\cdots,n)$. Let
\begin{eqnarray*}
f(z)&=&(p_1(z_1,z_2,\cdots, z_n),  p_2(z_{2})+f_{2}(z_{k}),\cdots,p_{k-1}(z_{k-1})+f_{k-1}(z_{k}),\\
&& p_k(z_{k}) ,p_{k+1}(z_{k+1})+f_{k+1}(z_{k}), \cdots, p_{n}(z_{n})+f_{n}(z_{k})).
\end{eqnarray*}
 If for any $z=(z_{1},z_{2},\cdots,z_{n})\in D_{p}^n\backslash\{0\}$, we have
\begin{eqnarray*}
&(1)&\frac{\partial p_1}{\partial z_1}\cdot\prod\limits_{j=2}^n
p_j'(z_j)\neq 0 ,\quad |z_j p_j''(z_j)|\leq |p_j'(z_j)|, j=2,\cdots, n;\\
&(2)&\sum\limits_{l=1}^n\bigg|z_1\frac{\partial^2 p_1}{\partial z_1\partial z_l}\bigg|\leq \bigg|\frac{\partial p_1}{\partial z_1}\bigg|;\\
&(3)&p_1\bigg|\frac{\frac{\partial p_1}{\partial
z_j}\cdot\frac{p_j''(z_j)}{p_j'(z_j)}} {\frac{\partial p_1}{\partial
z_1}}\bigg|
+p_1\sum\limits_{l=1}^n\bigg|\frac{\frac{\partial^2
p_1}{\partial z_j\partial z_l}} {\frac{\partial p_1}{\partial
z_1}}\bigg| \le
p_{j}|z_j|^{p_j-2}\bigg(1-\bigg|\frac{z_jp_j''(z_j)}{p_j'(z_j)}\bigg|\bigg),\\
&&(j=2,\cdots,k-1,k+1,\cdots, n-1);\\
&(4)&\sum\limits_{j=2,j\neq k}^{n}\bigg|\frac{f_j''(z_k)}{p_j'(z_j)}\bigg|p_{j}
 +\sum\limits_{j=2,j\neq k}^{n}\bigg|\frac{f_j'(z_k)}{p_j'(z_j)}\bigg|\bigg|\frac{p_k''(z_k)}{p_k'(z_k)}\bigg|p_{j}
 + \sum\limits_{j=2,j\neq k}^{n}\bigg|\frac{\frac{f_j''(z_k)}{p_j'(z_j)}\frac{\partial p_1}{\partial z_j}}{\frac{\partial p_1}{\partial z_1}}\bigg|p_1\\
\end{eqnarray*}
\begin{eqnarray*}
&& + \sum\limits_{j=2,j\neq k}^{n}\bigg|\frac{\frac{f_j'(z_k)}{p_j'(z_j)}\frac{p_k''(z_k)}{p_k'(z_k)}\frac{\partial p_1}{\partial z_j}}{\frac{\partial p_1}{\partial
 z_1}}\bigg|p_{1}
+\sum\limits_{l=1}^n\bigg|\frac{\frac{\partial^2
p_1}{\partial z_l\partial z_k}}{\frac{\partial p_1}{\partial
z_1}}\bigg|p_1
+\bigg|\frac{\frac{p_k''(z_k)}{p_k'(z_k)}\frac{\partial
p_1}{\partial z_k}}{\frac{\partial p_1}{\partial z_1}}\bigg| p_1\\
&&\le \bigg(1-\bigg|\frac{z_k p_k''(z_k)}{p_k'(z_k)}\bigg|\bigg)
p_{k}|z_k|^{p_k-2},
\end{eqnarray*}
then $f\in K(D^n_p)$.

Finally, we give a partial answer to Problem II by verifying the following theorem.

{\bf Theorem 4}\quad Suppose that $p\geq 2$. Let
$$
f(z)=(z_{1}+a_{1}z^{2}_{1}+a'_{1}z^{2}_{2},a_{2}z^{2}_{1}+z_{2}+a'_{2}z^{2}_{2}),
$$
where $z=(z_{1},z_{2})$. If $(a_{1},a_{2}, a'_{1}, a'_{2})$ satisfies the following conditions:
\begin{eqnarray}
4|a_{1}|+2|a'_{2}|+8|a_{1}||a'_{2}|+2|a_{2}|+4|a_{1}||a_{2}|+8|a'_{1}||a_{2}|+4|a_{1}||a_{2}|&\leq& 1,\label{liu4}\\
2|a_{1}|+2|a'_{1}|+4|a'_{1}||a'_{2}|+4|a'_{2}| +8|a_{1}||a'_{2}|+4|a'_{1}||a'_{2}|+8|a'_{1}||a_{2}|&\leq& 1.\label{liu5}
\end{eqnarray}
Then $f\in K(B^{2}_{p})$.

{\bf Proof}\quad By calculating the Fr\'{e}chet derivatives of $f(z)$ straightforwardly, we obtain
\begin{eqnarray*}
&&Df(z)=\begin{pmatrix}1+2a_{1}z_{1}&2a'_{1}z_{2}\\2a_{2}z_{1}&1+2a'_{2}z_{2}
\end{pmatrix},\quad
Df(z)^{-1}=\begin{pmatrix} \frac{1+2a'_{2}z_{2}}{A}
&\frac{-2a'_{1}z_{2}}{A}\\
\frac{-2a_{2}z_{1}}{A}
&\frac{1+2a_{1}z_{1}}{A}
\end{pmatrix},\\
&&D^{2}f(z)(b,b)=\begin{pmatrix}2a_{1}b_{1}&2a'_{1}b_{2}\\2a_{2}b_{1}&2a'_{2}b_{2}\\\end{pmatrix}
\begin{pmatrix}
b_{1}\\
b_{2}
\end{pmatrix}
=\begin{pmatrix}
2a_{1}b^{2}_{1}+2a'_{1}b^{2}_{2}\\
2a_{2}b^{2}_{1}+2a'_{2}b^{2}_{2}
\end{pmatrix},\\
&&Df(z)^{-1}D^{2}f(z)(b,b)=\begin{pmatrix}\frac{(1+2a'_{2}z_{2})(2a_{1}b^{2}_{1}+2a'_{1}b^{2}_{2})}
{A}
-\frac{2a'_{1}z_{2}(2a_{2}b^{2}_{1}+2a'_{2}b^{2}_{2})}{A}\\
-\frac{2a_{2}z_{1}(2a_{1}b^{2}_{1}+2a'_{1}b^{2}_{2})}{A}
+\frac{(1+2a_{1}z_{1})(2a_{2}b^{2}_{1}+2a'_{2}b^{2}_{2})}{A}
\end{pmatrix},
\end{eqnarray*}
where
\begin{eqnarray}
A=4a_{1}a'_{2}z_{1}z_{2}-4a'_{1}a_{2}z_{1}z_{2}+1+2a'_{2}z_{2}+2a_{1}z_{1}.\label{liu3}
\end{eqnarray}

Taking $z=(z_{1}, z_{2})\in B^2_p\setminus\{0\}, b=(b_{1}
,b_{2})\in C^2$ such that $\mbox{Re}\langle b, \frac{\partial
u}{\partial\overline{z}}\rangle =0$,  by the hypothesis of Theorem 2, we have
\begin{eqnarray*}
&&\frac{2}{p}J_{f}(z,b)\geq \he{j=1}{2}|b_{j}|^{2}|z_{j}|^{p-2}-\frac{2}{p}\mbox{Re}\langle Df(z)^{-1}D^{2}f(z)(b,b), \frac{\partial u}{\partial\overline{z}}\rangle\\
&=&\he{j=1}{2}|b_{j}|^{2}|z_{j}|^{p-2}-\mbox{Re}\bigg\{ \bigg(\frac{(1+2a'_{2}z_{2})(2a_{1}b^{2}_{1}+2a'_{1}b^{2}_{2})-2a'_{1}z_{2}(2a_{2}b^{2}_{1}+2a'_{2}b^{2}_{2})}{A} \bigg)\frac{|z_{1}|^{p}}{z_{1}}\\
&&+\bigg(\frac{(1+2a_{1}z_{1})(2a_{2}b^{2}_{1}+2a'_{2}b^{2}_{2})-2a_{2}z_{1}(2a_{1}b^{2}_{1}+2a'_{1}b^{2}_{2})}{A}\bigg)
\frac{|z_{2}|^{p}}{z_{2}}\bigg\}\\
&\geq &\he{j=1}{2}|b_{j}|^{2}|z_{j}|^{p-2}-\bigg\{ \bigg(\bigg|\frac{2a_{1}+4a_{1}a'_{2}z_{2}}{A}\bigg||b_{1}|^{2}
+\bigg|\frac{2a'_{1}+4a'_{1}a'_{2}z_{2}}{A}\bigg||b_{2}|^{2}+\bigg|\frac{4a'_{1}a_{2}z_{2}}{A}\bigg||b_{1}|^{2}\\
&&+\bigg|\frac{4a'_{1}a'_{2}}{A}\bigg||b_{2}|^{2}\bigg)|z_{1}|^{p-1}+\bigg(\bigg|\frac{4a_{1}a_{2}z_{1}}{A}\bigg||b_{1}|^{2}+\bigg|\frac{4a'_{1}a_{2}z_{1}}{A}\bigg||b_{2}|^{2}\\
\end{eqnarray*}
\begin{eqnarray*}
&&+\bigg|\frac{2a_{2}+4a_{1}a_{2}z_{1}}{A}\bigg||b_{1}|^{2}
+\bigg|\frac{2a'_{2}+4a_{1}a'_{2}z_{1}}{A}\bigg||b_{2}|^{2}\bigg)|z_{2}|^{p-1}\bigg\}\\
&\geq & |b_{1}|^{2}|z_{1}|^{p-2}\bigg(1-\frac{2|a_{1}|+4|a_{1}||a'_{2}|+2|a_{2}|+4|a_{1}||a_{2}|+4|a'_{1}||a_{2}|+4|a_{1}||a_{2}|}
{1-4|a_{1}||a'_{2}|-4|a'_{1}||a_{2}|-2|a'_{2}|-2|a_{1}|}\bigg)\\
&&+|b_{2}|^{2}|z_{2}|^{p-2}\bigg(1-\frac{2|a'_{1}|+4|a'_{1}||a'_{2}|+2|a'_{2}| +4|a_{1}||a'_{2}|+4|a'_{1}||a'_{2}|+4|a'_{1}||a_{2}|}{1-4|a_{1}||a'_{2}|-4|a'_{1}||a_{2}|-2|a'_{2}|-2|a_{1}|}\bigg)\\
&\geq &0,
\end{eqnarray*}
hence it follows from Theorem A that $f(z)\in K(B^{2}_{p})$. \hfill $\Box$

\end{document}